\documentclass[10pt]{article}
\usepackage[english,activeacute]{babel}
\usepackage{epsfig}

\usepackage[latin1]{inputenc}
\usepackage[T1]{fontenc}
\usepackage{stmaryrd}
\usepackage{amsmath,amsfonts,amssymb,mathrsfs,amsthm}
\usepackage{xcolor}
\usepackage{dsfont}
\usepackage{graphicx}
\usepackage[mathscr]{eucal}
\usepackage{makeidx}
\usepackage{verbatim}
\usepackage{graphics,graphicx}
\usepackage{textcomp}
\usepackage{float}
\usepackage{tikz}
\usepackage[colorlinks=true, linkcolor=purple, citecolor=purple]{hyperref}

\newtheorem{lemma}{Lemma}[section]
\newtheorem{theorem}[lemma]{Theorem}
\newtheorem{proposition}[lemma]{Proposition}
\newtheorem{corollary}[lemma]{Corollary}
\newtheorem{assumption}[lemma]{Assumption}

\theoremstyle{definition}

\newtheorem{definition}[lemma]{Definition}
\newtheorem{remark}[lemma]{Remark}

\makeatletter
\def\keywords{
    \vspace{1ex}
    \noindent
    \if@twocolumn
      \small{\bf  Keywords}\/---$\!$    \else
      \begin{center}\small\ {\bf Keywords}\end{center}\quotation\small
    \fi}
\def\endkeywords{\vspace{0.6em}\par\if@twocolumn\else\endquotation\fi
    \normalsize\rm}
\makeatother

\renewcommand{\O}{\ensuremath{\mathcal O}}

\renewcommand{\L}{\ensuremath{\mathcal L}}

\newcommand{\calS}{\ensuremath{\mathcal S}}

\newcommand{\calC}{\ensuremath{\mathcal C}}

\newcommand{\x}{\ensuremath{\bf x}}

\DeclareMathOperator{\Sp}{Sp}

\DeclareMathOperator{\comp}{comp}

\newcommand{\mb}[1]{\ensuremath{\mathbb{#1}}}
\newcommand{\N}{{\mb{N}}}

\newcommand{\R}{{\mb{R}}}
\newcommand{\C}{{\mb{C}}}

\newcommand{\y}{\ensuremath{y}}

\newcommand{\eps}{\varepsilon}

\newcommand{\M}{\ensuremath{\mathcal M}}

\let \Re \relax
\DeclareMathOperator{\Re}{Re}

\newcommand{\ovl}[1]{\overline{#1}}
\newcommand{\udl}[1]{\underline{#1}}

\DeclareMathOperator{\supp}{supp}

\DeclareMathOperator{\Op}{Op}

\newcommand{\transp}{\ensuremath{\phantom{}^{t}}}

\renewcommand{\d}{\ensuremath{\partial}}

\def\x{x}

\newcommand{\E}{\mathscr E}

\newcommand{\rouge}[1]{{\color{red} #1}}

\newcommand\bna{\begin{eqnarray*}}
\newcommand\ena{\end{eqnarray*}}

\def\bal#1\nal{\begin{align*}#1\end{align*}}
\def\baln#1\naln{\begin{align}#1\end{align}}

\newcommand\bnan{\begin{eqnarray}}
\newcommand\enan{\end{eqnarray}}

\newcommand\bnp{\begin{proof}}
\newcommand\enp{\end{proof}}

\newcommand\bneq{\begin{eqnarray*}\left\lbrace \begin{array}{rcl}}
\newcommand\eneq{\end{array} \right.\end{eqnarray*}}
\newcommand\bneqn{\begin{eqnarray}\left\lbrace \begin{array}{rcl}}
\newcommand\eneqn{\end{array} \right.\end{eqnarray}}

 \numberwithin{equation}{section} 
 
\newcommand\grando[1]{\mathcal{O}\left(#1\right)}

\renewcommand\M{[0,L]}

\newcommand\nor[2]{\left\|#1\right\|_{#2}}

\newcommand\trans{{}^t \!}
\newcommand\petito[1]{\mathbf{o}(#1)}
\newcommand\e{\varepsilon}

\newcommand{\qe}{q_{\eps}}
\newcommand{\qen}{q_{\eps_n}}

\newcommand\V{\mathcal{V}}

\newcommand{\xzero}{\ensuremath{\mathbf{x}_0}}

\begin{document}
\title{
Uniform observation of semiclassical Schr\"odinger eigenfunctions on an interval
}

\author{Camille Laurent\footnote{CNRS UMR 7598 and Sorbonne Universit\'es UPMC Univ Paris 06, Laboratoire Jacques-Louis Lions, F-75005, Paris, France, email: camille.laurent@sorbonne-universite.fr} and Matthieu L\'eautaud\footnote{Laboratoire de Math\'ematiques d'Orsay, UMR 8628, Universit\'e Paris-Saclay, CNRS, B\^atiment 307, 91405 Orsay Cedex France, email: matthieu.leautaud@universite-paris-saclay.fr.}
}

\date{}

\maketitle

\begin{abstract}
We consider eigenfunctions of a semiclassical Schr\"odinger operator on an interval, with a single-well type potential and Dirichlet boundary conditions. We give upper/lower bounds on the $L^2$ density of the eigenfunctions that are uniform in both semiclassical and high energy limits. These bounds are optimal and are used in an essential way in the companion paper~\cite{LL:20-1D} in application to a controllability problem. The proofs rely on Agmon estimates and a Gronwall type argument in the classically forbidden region, and on the description of semiclassical measures for boundary value problems in the classically allowed region.
Limited regularity for the potential is assumed.
\end{abstract}


\tableofcontents

\section{Introduction and main results}
\begin{figure}[h!]
  \begin{center}
\begin{tikzpicture}
\draw[->] (-5,0) -- (6,0) node[right] {$x$}; 
\draw[->] (-3,-1) -- (-3,5) ; 
\draw[-] (3,-1) -- (3,5);  

\begin{scope}
\clip (-5,-1) rectangle (4,5);
\draw[color=red,samples=100] plot ({\x},{1/5*(\x-1)*(\x-1)+1});
\end{scope}

\draw[color=blue,dashed] (1,-0.5) -- (1,1.5);
\draw[color=blue,dashed] (-3.5,1) -- (3.5,1);
\node[below left=0.1cm] at (-3,-0) {$0\strut$};
\node[below left=0.1cm] at (3,0) {$L\strut$};
\node[color=red,below right=0.1cm] at (3.8,3) {$V(x)\strut$};
\node[color=blue,above left=0.1cm] at (-3,1) {$E_0=V(\xzero)\strut$};
\node[color=blue,below left=0.1cm] at (1,0) {$\xzero\strut$};

\end{tikzpicture}
    \caption{A typical potential $V$ satisfying Assumption~\ref{assumptions} (and thus Assumption~\ref{assumptions-0})}
    \label{f:geom-setting}
 \end{center}
\end{figure}
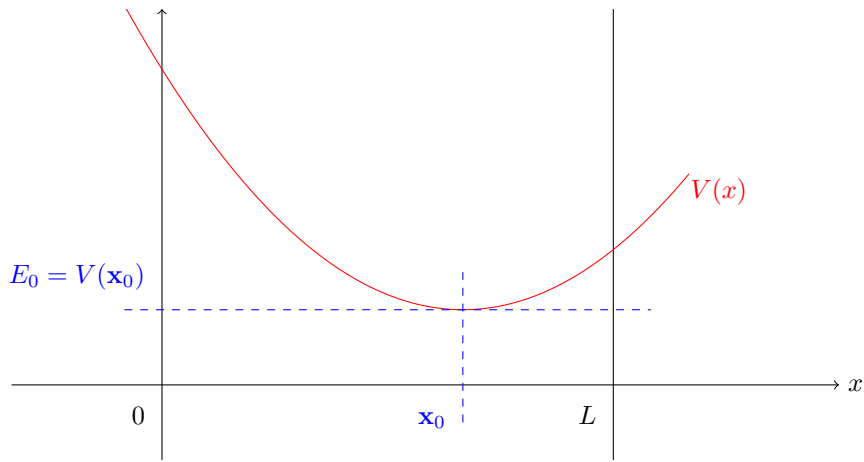

We investigate the localization of eigenfunctions of the semiclassical Schr\"odinger operator
 \begin{equation}
 \label{e:def-Peps}
 P_\eps := - \eps^2 \d_x^2 +V_\eps(x) ,
 \end{equation}
 on the interval $[0,L]$, with Dirichlet boundary conditions, where $V_\eps : [0,L]\to \R$ is a family of real-valued bounded potentials.  
In this setting, for any $\eps>0$, the operator $P_\eps$ endowed with domain $D(P_\eps)= H^2( [0,L]) \cap H^1_0( [0,L])$ is a selfadjoint operator on $L^2(0,L)$, with compact resolvents. Its spectrum $\Sp(P_\eps)$ thus consists only in countably many real eigenvalues with finite multiplicity (equal to $1$ since this is a $1D$ problem).
 We are concerned with properties of eigenfunctions of $P_\eps$, that is to say, solutions $\psi$ to  
 \begin{align}
\label{e:eignfct-E-bis}
P_\eps \psi = E \psi , \qquad
 \psi \in H^2( [0,L]) \cap H^1_0( [0,L]), \qquad  \nor{\psi}{L^2(\M)} =1 ,
\end{align}
where, as already mentioned, $E$ is necessary a real number (depending on $\eps$).
We shall further assume that the potentials $V_\eps$ converge to a fixed potential $V$. The assumptions we make on $V_\eps$ and $V$ are one of the two following.
  \begin{assumption}
\label{assumptions-0} 
Assume
\begin{itemize}
\item $V \in C^0([0,L];\R)$, $V_\eps \in L^\infty(0,L;\R)$ are real valued and $ \|V - V_\eps\|_{L^\infty(0,L)}\to 0$; 
\item there is $x= \xzero \in (0,L)$ such that $V$ is strictly decreasing on $[0,\xzero]$ and strictly increasing on $[\xzero,L]$.
 \end{itemize}
\end{assumption}
  \begin{assumption}
\label{assumptions} 
Assume
\begin{itemize}
\item $V_\eps,V \in C^1([0,L];\R)$ are real valued and   $ \|V - V_\eps\|_{C^1([0,L])}\to 0$; 
\item the only $x \in [0,L]$ such that $V'(x) = 0$ is $x= \xzero \in (0,L)$ and $V(\xzero) = \min_{[0,L]}V$.
 \end{itemize}
\end{assumption}
The typical shape of the potential $V$ is illustrated on Figure~\ref{f:geom-setting}.

Note that Assumption~\ref{assumptions}  implies Assumption~\ref{assumptions-0}. 
We shall also write alternatively $V_\eps(x) = V(x)+  \qe (x)$ with $\qe \to 0$ in $L^\infty$ or $C^1$ topology as $\eps \to 0$.
That is to say, we consider the single well problem on the interval.  
We denote by $E_0$ the ground state energy, that is to say
$$
E_0 =\min_{x\in \M} (V(x)) = V(\xzero) . 
$$
The classically allowed region at energy $E$ for the potential $V$ is defined by: 
$$
 K_E = \{x \in \M , V(x) \leq E \} ,
$$
and the Agmon distance (see e.g.~\cite[Chapter~3]{Helffer:booksemiclassic}) to the set $K_E$ at the energy level $E$ by
\begin{align}
\label{defAgmonbis}
d_{A,E}(x) : = \inf_{y \in K_E} \left|\int_{y}^{x} \sqrt{\left(V(s) - E \right)_+ } ds\right|=  \left|\int_{y_E}^{x} \sqrt{\left(V(s) - E \right)_+ } ds\right| , \quad \text{ if } E\geq E_0,
\end{align}
 where $\left(V(x) - E \right)_+ = \max \left(V(x) - E , 0\right)$ and where $y_E$ is any point in $K_E$.
Note in particular that $d_{A,E}$ vanishes identically on $K_E$ (and only on this set).
If $E <E_0$, we have $K_E=\emptyset$ so that the Agmon distance above is not well-defined; in that case, we shall use the convention that 
$$
d_{A,E}(x) =d_{A,E_0}(x) , \quad \text{ if } E\leq E_0 .
$$
This is the appropriate convention since, 
 if $\psi$ and $E\in \R$ satisfy~\eqref{e:eignfct-E-bis}, the $L^2$ inner product of~\eqref{e:eignfct-E-bis} with $\psi$ yields 
\begin{align}
\label{e:E=nor-psi}
E = \eps^2\nor{\psi'}{L^2(\M)}^2 + \int_{\M} (V+\qe) |\psi|^2 ,
\end{align}
and thus, under Assumption~\ref{assumptions-0},
\begin{align}
\label{e:E=nor-psi-bis}
E \in \Sp(P_\eps) \implies E\geq E_0 - \|\qe\|_\infty \to_{\eps\to 0^+} E_0 .
\end{align}

Under Assumption~\ref{assumptions}, we prove upper and lower bounds that, roughly speaking, say that solutions of $P_\eps \psi = E \psi$ behave, in the sense of $L^2$-density, like $|\psi(x)|\sim e^{-\frac{d_{A,E}(x)}{\eps}}$ up to some loss $e^{\frac{\delta}{\eps}}$.
The upper bounds on the eigenfunctions of $P_{\e}$ are expressed under the form of uniform Agmon estimates. 

\begin{theorem}[Upper bounds on eigenfunctions: uniform Agmon type estimates]
\label{t:uniform-agmon}
Let $V,V_\eps$  satisfy Assumption~\ref{assumptions-0}.
Then, for all $\delta>0$ there exist $\eps_0 = \eps_0(\delta) \in (0,1]$ such that for all $E \in \R$ and $\psi$ solution to~\eqref{e:eignfct-E-bis},
 we have for all $\eps< \eps_0$
\begin{align}
\label{e:agmon-H1}
 \nor{ e^{\frac{d_{A,E}}{\eps}}\frac{\eps}{\sqrt{|E|+1}}\psi'}{L^2} + \nor{ e^{\frac{d_{A,E}}{\eps}}\psi}{L^2} \leq e^{\frac{\delta}{\eps}}  .\\
\label{e:Agmon-boundary}
\frac{\eps}{\sqrt{|E|+1}}|\psi'(0)| \leq e^{- \frac{d_{A,E}(0)-\delta}{\eps}} , \quad \frac{\eps}{\sqrt{|E|+1}}|\psi'(L)| \leq e^{- \frac{d_{A,E}(L)-\delta}{\eps}}  .
\end{align}
\end{theorem}

The main result of this note is the following converse estimate.

\begin{theorem}[Lower bounds on eigenfunctions]
\label{t:allibert-lower-uniform}
Let $V,V_\eps$ satisfy Assumption~\ref{assumptions}.
Then, for any interval $U\subset [0,L]$ with nonempty interior and any $\delta >0$, there is $\eps_0>0$ such that for all 
$E \in \R$ and $\psi$ solution to~\eqref{e:eignfct-E-bis},
 we have for all $\eps<\eps_0$, 
\begin{align}
&\nor{\psi}{L^2(U)} \geq e^{-\frac{1}{\eps}(d_{A,E}(U) + \delta)} ,\quad d_{A,E}(U) = \inf_{x\in U}d_{A,E}(x), 
\label{e:internal-obs} \\
&\frac{\eps}{\sqrt{|E|+1}} |\psi'(0)|  \geq e^{-\frac{1}{\eps}(d_{A,E}(0) + \delta)}, \quad \frac{\eps}{\sqrt{|E|+1}} |\psi'(L)|  \geq e^{-\frac{1}{\eps}(d_{A,E}(L) + \delta)} \label{e:bound-obs}  .
\end{align}
\end{theorem}
Note that this lower bound is as precise as the upper bound~\eqref{e:agmon-H1} (except for the $\delta$ loss) and thus essentially optimal. Also, in these estimates, the loss $e^{-\frac{\delta}{\eps}}$ can be removed/improved in several situations (see e.g. Proposition~\ref{l:geom-control} in the classically allowed region).

Remark that Theorems~\ref{t:allibert-lower-uniform} and~\ref{t:uniform-agmon} are counterparts one to the other. They state essentially that, in this very particular one dimensional setting, an eigenfunction $\psi$ associated to the energy $E$ satisfies $|\psi(x)|\sim e^{-\frac{d_{A,E}(x)}{\eps}}$ in the sense of $L^2$ density (and that this is uniform in $E, x,\eps$).

Notice finally that, under Assumption~\ref{assumptions-0}, the set $K_E$ is an interval given for $E\geq E_0$ by $K_E=[x_-(E),x_+(E)] \subset [0,L]$, where $x_\pm(E)$ are defined precisely below.
\begin{definition}
\label{d:def-xpm}
For $E\geq E_0$, set 
\begin{itemize}
\item  $x_-(E)$ the solution to $V(x_-(E))=E$ which is $\leq \xzero$ for $E \leq V(0)$, and  $x_-(E)=0$ for $E \geq V(0)$,   
\item $x_+(E)$ the solution to $V(x_+(E))=E$ which is $\geq \xzero$ for $E \leq V(L)$, and  $x_+(E)=L$ for $E \geq V(L)$,
\end{itemize}(with $\xzero=x_-(E_0)=x_+(E_0)$ if $E=E_{0}$).
\end{definition}

\bigskip
The proof of Theorem~\ref{t:allibert-lower-uniform} relies on an explicit expression of semiclassical measures in the present context, which is of its own interest.

\begin{theorem}
\label{thmmeasurex}
Assume that $V_\eps, V$ satisfy Assumption~\ref{assumptions}.
Suppose that $\eps_n\to 0$, $E_n\to E_* \in \R\cup \{+\infty\}$ as $n \to +\infty$, and $\psi_n$ solves
 \begin{align}
\label{e:eignfct-E-n}
(P_{\eps_n}  - E_n) \psi_n  = r_n, \qquad
 \psi_n \in H^2( [0,L]) \cap H^1_0( [0,L]), \qquad  \nor{\psi_n}{L^2(\M)} =1 ,
\end{align}
where $\|r_n\|_{L^2(0,L)} = \petito{\eps_n}$.
 Then, in the sense of weak$-*$ convergence of measures, we have $|\psi_n(x)|^2 dx \rightharpoonup \mathfrak{m}_{E_*}$ for a nonnegative Radon measure $\mathfrak{m}_{E_*}$ on $[0,L]$ explicitely given by
\begin{align*}
&\mathfrak{m}_{E_*}= C_{E_*} \frac{\mathds{1}_{(x_-(E) , x_+(E))}(x)dx }{\sqrt{(E-V(x))_+}}  , \quad \text{ if }E_0<E_*<+\infty , \quad\text{ with } C_{E_*} = \left(  \int_{x_-(E_{*})}^{x_+(E_{*})}\frac{dx }{\sqrt{E_{*}-V(x)}} \right)^{-1} , \\
&\mathfrak{m}_{E_*}= \delta_{\xzero} , \quad \text{ if }E_*=E_0 , \\
&\mathfrak{m}_{E_*}= \frac{dx}{L}, \quad \text{ if } E_*=+\infty .
\end{align*}
Moreover, in $\R$ we have 
\begin{align*}
&|\eps_n\psi_n'(0)|^2 \to  2 C_{E_*}\sqrt{E_*-V(0)} \mathds{1}_{V(0) <E_*} , \quad |\eps_n\psi_n'(L)|^2 \to  2 C_{E_*}\sqrt{E_*-V(L)} \mathds{1}_{V(L) <E_*} , \quad \text{ if }E_* < +\infty , \\
&E_{n}^{-1}|\eps_n\psi_n'(0)|^2 \to  \frac{2}{L},\quad E_{n}^{-1}|\eps_n\psi_n'(L)|^2 \to  \frac{2}{L}, \quad \text{ if } E_*=+\infty .
\end{align*}
\end{theorem}
Several remarks are in order.
First, for a given $E_*$, the uniqueness of the limit measure implies that the whole sequence $|\psi_n(x)|^2 dx$ converges. This is an extremely rare situation (probably linked to the simplicity of the spectrum and the regularity of the spectral gap in this 1D situation, but we do not use this information here).

Second, this theorem only describes the limit measures of $|\psi_n(x)|^2 dx$. The latter are  projections on the $x$-space of the semiclassical measure that live in the phase-space $(x,\xi) \in [0,L]\times \R$, and are as well described explicitly in the proof of Theorem~\ref{thmmeasurex}. Their expression is slightly less readable, so that we decided not to write them here.

Other possible approaches to this problem (which could in principle also lead to statements like those of Theorems~\ref{t:allibert-lower-uniform} and~\ref{thmmeasurex}) include WKB expansions (at least to leading order), see e.g.~\cite[pp139-143]{GS:94} for the single well problem in $\R$ or~\cite{Dui:74} (in a much more general seeting), or ODE methods see e.g. \cite[Section~6~pp~190--198]{Olver:book},~\cite[Theorems~4.5 and~4.6]{BS:book} or~\cite{FF:02}.

\bigskip
The study of eigenvalues and eigenfunctions for $1D$ Schr\"odinger operators in the semiclassical limit is a classical topics; we refer e.g. to the seminal papers of Simon~\cite{Simon:83} and Helffer-Sj\"ostrand~\cite{HS:84} for the bottom energy and Helffer-Robert~\cite{HR:84} for higher energies, as well as the books by Helffer and Dimassi-Sj\"ostrand~\cite{Helffer:booksemiclassic,DS:book}. 
In particular, the proof of Theorem~\ref{t:uniform-agmon} consists in a rather classical Agmon estimates~\cite{HS:84,Helffer:booksemiclassic,DS:book}, and we essentially need to check here the limited regularity of the potential and the uniform dependence on the energy levels $E$. This uniformity is necessary for the proof of Theorem~1.6 in~\cite{LL:20-1D}.

The literature on lower bounds (such as given in Theorem~\ref{t:allibert-lower-uniform}) and semiclassical measures (such as given in Theorem~\ref{thmmeasurex}) for a boundary value problem is slightly poorer. We mention the article~\cite{Allibert:98} where an analogue of Theorem~\ref{t:allibert-lower-uniform} is stated in which the lower bounds in the right hand-sides of~\eqref{e:internal-obs} and~\eqref{e:bound-obs} is given in terms of the Agmon distance to the ground energy $d_{A,E_0}$.
Similar (but less precise) estimates have been also used by the authors in~\cite{LL:18,LL:18vanish} for applications to eigenfunctions on surfaces of revolution.

The exponential bounds obtained in both Theorem~\ref{t:uniform-agmon} and~\ref{t:allibert-lower-uniform} could certainly be refined under additional assumptions (analyticity of $V_\eps=V$, non degeneracy of $V$ at $\xzero$...), especially for the bottom energy $E_0$, using e.g. some of the techniques developed in~\cite{HS:84,HS:86,Helffer:booksemiclassic,DS:book,HN:06}.

Note finally that there are very few situations in which semiclassical measures of eigenfunctions/quasimodes can be described explicitely; see e.g.~\cite{JakobsonTori97} on the torus or~\cite{ALM:cras} on the disk. It is therefore satisfactory to be able to express all semiclassical measures in this very simple geometric situation.
We refer to~\cite[Section~4]{HMR:87} (relying on~\cite{Dui:74}) for a related statement in a boundaryless setting with $V_\eps=V$ smooth, linked to quantum ergodicity. Note by the way that the proof of Theorem~\ref{thmmeasurex} below implies in particular that the operator~\eqref{e:def-Peps} is quantum unique ergodic at all energy levels under Assumption~\ref{assumptions}.

\bigskip
The plan of the article is thus as follows. Section~\ref{s:proofs} is devoted to the proofs of the above results.
The proof of Theorem~\ref{t:uniform-agmon}, consequence of Agmon estimates, is first given in Section~\ref{s:proofAgmon} below as a warmup. 
Then, we focus on the proof of Theorem~\ref{t:allibert-lower-uniform}, which relies on three key lemmata:
\begin{itemize}
\item a geometric control estimate in the classically allowed region, proved in Section~\ref{s:classically-allowed-section}. The latter essentially reduces to the description of semiclassical measures as stated in Theorem~\ref{thmmeasurex}, and Section~\ref{s:classically-allowed-section} is thus dedicated to the proof of Theorem~\ref{thmmeasurex};
\item a tunneling estimate into the classically forbidden region (inspired by~\cite{Allibert:98}), with sharp tunneling rate, proved in Section~\ref{s:class-forbid};
\item a rough Gronwall estimate used to patch the previous two estimates in the transition between the classically allowed and forbidden regions (that is, near the two turning points), also proved in Section~\ref{s:class-forbid}.
\end{itemize} 
The last two points use arguments inspired by~\cite[Section~3.2 pp1541-1546]{Allibert:98}. There are three main differences with that reference. First, we have $d_{A,E}$ in the exponent of Theorem~\ref{t:allibert-lower-uniform}, where Allibert only had $d_{A,E_0}$. Second, our estimate is uniform with respect to the energy level $E$. Third, the potential has limited regularity and can be perturbed by lower order terms (denoted $\qe$ here).
This uniformity is actually a source of some complications in the proofs. Yet, it is necessary for the proof of the cost of controllability in Theorem~1.6 in~\cite{LL:20-1D}.
We finally prove Theorem~\ref{t:allibert-lower-uniform} from the three key lemmata in Section~\ref{l:end-of-proof}.

\medskip
Section~\ref{s:app-msc} is devoted to the proof of several technical properties of semiclassical measures for boundary-value problems (and in dimension one only), that are prerequisites to the proof of Theorem~\ref{thmmeasurex}. The results are summarized in Proposition \ref{propmeasureapp}.

The plan of Section~\ref{s:app-msc} is as follows.
We start by proving a priori estimate and the so-called hidden regularity of traces in Section~\ref{s:reg-of-traces}. This allows to define semiclassical measures associated to the eigenfunctions $\psi_n(x)$ (as well as limits of the Neumann traces), that are lifts to the phase space $(x,\xi)\in [0,L]\times \R$ of the measures $\mathfrak{m}_{E_*}$ appearing in Theorem~\ref{thmmeasurex}. 
We then prove that these semiclassical measures are supported on the energy layer $\{\xi^2+V(x)=E_*\}$ in Section~\ref{s:loc-char-set}. 
Next, we prove in Section~\ref{s:propagation} that the measure satisfies an appropriate transport equation (charged at the boundary).
Invariance properties near the boundary are finally deduced in Section~\ref{s:invar-boundar}.

Most arguments in Section~\ref{s:app-msc} are essentially inspired from the seminal paper of G\'erard and Leichtnam~\cite{GL:93}, where eigenfunctions of the Laplace operator are considered in any dimension, in domains with boundary having limited smoothness. We believe it is useful to provide here with a detailed argument in our context for two reasons. First, the results of \cite{GL:93} do not apply here since they only deal with the flat Laplacian without potential. Second, the proofs of \cite{GL:93} (as well as other references on boundary propagation for semiclassical measures, e.g.~\cite{Leb:96,Burq:97,Burq:97b,RZ:09}) are highly technical because of the geometry and the weak regularity of the boundary. Many arguments simplify considerably in our 1D context. We thus take this as an opportunity to write a proof as detailed and pedagogical as possible, which we hope can be read as an elementary introduction to boundary propagation.

Note that although the problem is one dimensional, the fact that we consider a semiclassical Schr\"odinger operator makes it a very good toy model that encompasses part of the richness of propagation theory for boundary value problems~\cite{MS:78}. Indeed, we shall see that elliptic, hyperbolic and glancing points all arise on the energy layer $\xi^2+V(x)=E_*$ for certain values of the energy $E_*$ (see Section~\ref{s:invar-boundar}).

 Note finally that all proofs of the present article are completely self-contained except for the standard semiclassical calculus in $\R$. 

\bigskip
\noindent
{\em Acknowledgements.} 
The second author is partially supported by the Agence Nationale de la Recherche under grants SALVE ANR-19-CE40-0004 and ADYCT (ANR-20-CE40-0017).

The authors would like to thank Bernard Helffer for his comments on a preliminary version of this work and for and pointing out several references, and the anonymous referee for her/his useful remarks, which helped to improve the quality of the article. 

\section{Proofs}
\label{s:proofs}
Before turning to the proofs, we start with two simple remarks that will be used along the proofs. The first remark aims at reducing the proofs to the energies $E$ that are $\geq E_0$.

The first remark concerns the a priori regularity of the functions $x_\pm$ of Definition~\ref{d:def-xpm} and $d_{A,E}$ defined in~\eqref{defAgmonbis}.
\begin{lemma}
\label{l:regularity-x-da}
Under Assumption~\ref{assumptions-0}, the functions $x_\pm : [E_0, \infty) \to [0,L]$ are uniformly continuous function.
The function $\R \times [0,L] \to \R$ defined by $(E,x) \mapsto d_{A,E}(x)$ is uniformly continuous and $x \mapsto d_{A,E}(x)$ is $C$-Lipschitz with $C$ independent of $E$.
\end{lemma}
\bnp
The first statement comes from continuity of $V^{-1}$ on the compact $[\xzero, L]$ (and similarly on $[0,\xzero]$). The second statement follows from the explicit expression
\begin{equation*}
\begin{array}{ll}
 d_{A,E}(x) = \int_{x_+(E)}^{x} \sqrt{V(s) - E} ds ,& \quad \text{if} \quad E\geq E_0 , x \geq x_+(E) , \\
d_{A,E}(x) = 0 ,& \quad \text{if} \quad E\geq E_0 , x \in [ x_- (E),x_+(E)] , \\
 d_{A,E}(x) = \int_{x}^{x_-(E)} \sqrt{V(s) - E} ds ,& \quad \text{if} \quad E\geq E_0 , x \leq x_-(E) , \\
  d_{A,E}(x) =  d_{A,E_0}(x)  =\left|\int_{\xzero}^{x} \sqrt{V(s) - E} ds\right|  ,& \quad \text{if} \quad E\leq  E_0 , x \in [0,L] ,
\end{array}
\end{equation*}
and in particular, $d_{A,E}(x)=0$ for $E \geq \max V$ and $d_{A,E}(x) =  d_{A,E_0}(x)$ for $E\leq  E_0$. 
Moreover, we see that $d_{A,E}$ is $C$-Lipschitz with $C= \max\left\{ \sqrt{V(x) - E} , E\in [E_0, \max V], x \in [0,L]\right\}$.
\enp

The second remark concerns the reduction of the statements for all energy levels $E\in \R$ to only $E\geq E_0$.
\begin{remark}
\label{r:EgeqE0}
We notice that it suffices to prove the statements of Theorems~\ref{t:uniform-agmon}--\ref{t:allibert-lower-uniform} for $E\geq E_0$ (and not for all $E\in \R$). 

Indeed, if $P_\eps \psi=E \psi$, and if we set $E_\eps= E+\|\qe\|_\infty$, we then have $E_\eps\geq E_0$ from~\eqref{e:E=nor-psi-bis}. 
Moreover, with $\tilde{P}_\eps =P_\eps  +\|\qe\|_\infty$ (which is equal to $P_\eps$ with $\qe$ replaced by $\tilde \qe = \qe + \|\qe\|_\infty \geq 0$ which is such that $\|\tilde \qe\|_{L^\infty}\to 0$ under Assumption~\ref{assumptions-0} or $\|\tilde \qe\|_{C^1}\to 0$ under Assumption~\ref{assumptions}) we have $\tilde{P}_\eps \psi = E_\eps \psi$. 
 
The results of Theorems~\ref{t:uniform-agmon}--\ref{t:allibert-lower-uniform} apply to $\tilde{P}_\eps$ and $E_\eps\geq E_0$ with $d_{A,E}$ replaced by $d_{A,E_\eps}$.
The conclusion for all $E\in \R$ follows from Lemma~\ref{l:regularity-x-da} below:
 for any $\delta>0$ there is $\eps_0>0$ such that for all $d_{A,E}-\delta \leq d_{A,E_\eps} \leq d_{A,E} + \delta$ uniformly on $x \in [0,L]$ and $\eps <\eps_0$.
\end{remark}

\subsection{Uniform Agmon estimates: Proof of Theorem~\ref{t:uniform-agmon}}
\label{s:proofAgmon}

We follow e.g.~\cite{Helffer:booksemiclassic} for the proof of Theorem~\ref{t:uniform-agmon}.
  
\bnp[Proof of Theorem~\ref{t:uniform-agmon}]
Notice first that according to Remark~\ref{r:EgeqE0}, it suffices to consider $E\geq E_0$.
Next, consider the range $E \geq \max_{\M} V $. In that case, $(V-E)_+=0$ and the Agmon distance $d_{A,E}$ vanishes identically on $[0,L]$. Hence the statement~\eqref{e:agmon-H1} writes $$\frac{\eps}{\sqrt{|E|+1}}\nor{\psi'}{L^2} + \|\psi \|_{L^2} \leq e^{\frac{\delta}{\eps}},$$ which, for $\eps$ sufficiently small, is a consequence of $\nor{\psi}{L^2(\M)} =1$ together with 
$$
\eps^2\nor{\psi'}{L^2(\M)}^2  \leq (|E|+ \nor{V}{\infty} +1) ,
$$
which follows from~\eqref{e:E=nor-psi}.
Next, Estimate~\eqref{e:Agmon-boundary} holds uniformly on compact sets of energies $E$ as a consequence of the hidden regularity Estimate~\eqref{e:estim-bebete-trace} in Lemma~\ref{l:preliminary-trace} below, taken for $h=\eps$ and $\mathscr{V} =\mathscr{V}_1=V_\eps -E$.
For $E\geq 1$, we write Estimate~\eqref{e:estim-bebete-trace} for $h=\frac{\eps}{\sqrt{E}}$,$\mathscr{V}_{1} =\frac{V_\eps}{E} =\frac{h^2}{\eps^2}V_\eps$ and  $\mathscr{V}_{2}  =-1$. This  implies that $\frac{\eps}{\sqrt{E}}|\psi'(0)|  = h |\psi'(0)|  \leq C h^{-1}\|\mathscr{V}_{1}\|_{L^\infty}+C  \leq \frac{h}{\eps^2}\|V_\eps\|_{L^\infty}+C  \leq \eps^{-2}C_{V,\qe}$ uniformly in $E, \eps$, and in particular \eqref{e:Agmon-boundary} holds in this range of energies.

\medskip
We finally consider the most substantial case, namely $E \in [\min_{\M} V-1 ,\max_{\M} V ]$, and proceed with the proof of the Agmon estimates.
We start with the following integration by parts formula. 
For all $\phi \in W^{1,\infty}(0,L)$ and $u \in H^2 \cap H^1_0(0,L)$ we have 
$$
\int_0^L \left( \eps^2 |\d_x(e^{\phi/\eps}u)|^2 - |\d_x\phi|^2 e^{2\phi/\eps}|u|^2 \right)  = \Re \int_0^L e^{2\phi/\eps} (-\eps^2 \d_x^2u) \ovl{u}. 
$$
We use this identity with $u=\psi$ a solution to $-\eps^2 \psi'' + V\psi + \qe \psi =P_\eps \psi = E \psi$. This yields 
$$
\int_0^L  \eps^2 |\d_x(e^{\phi/\eps}\psi)|^2 +  \int_0^L (V -E - |\d_x\phi|^2 + \qe )e^{2\phi/\eps}|\psi|^2   = 0.
$$
We now write $(0,L) = \Omega^+_\alpha \sqcup \Omega^-_\alpha$ with $\Omega^+_\alpha = \{V-E \geq \alpha^2\}$ and $\Omega^-_\alpha= \{V-E < \alpha^2\}$ for some $0<\alpha\leq 1$ to be chosen later.  We obtain
\begin{align}
\label{e:agmon-prelim}
\int_0^L  \eps^2 |\d_x(e^{\phi/\eps}\psi)|^2 +  \int_{\Omega_\alpha^+} (V -E - |\d_x\phi|^2 + \qe )e^{2\phi/\eps}|\psi|^2 \leq \sup_{\Omega_\alpha^-} \left|V -E - |\d_x\phi|^2 + \qe \right|\int_{\Omega_\alpha^-}e^{2\phi/\eps}|\psi|^2 .
\end{align}
We now choose the weight $\phi = (1-\delta)d_{A,E}$ for $\delta \in (0,1)$ (where $d_{A,E}$ is defined in \eqref{defAgmonbis} and is Lipschitz continuous according to Lemma~\ref{l:regularity-x-da}). 

On $\Omega_\alpha^+$, noticing that $|d_{A,E}'|^2 = (V-E)_+= V-E$, we have 
$$
V -E - |\d_x\phi|^2 + \qe = (V -E)(1  -(1-\delta)^2) + \qe
 \geq \alpha^2\delta(2-\delta) - \nor{\qe}{\infty} ,
$$
hence providing with a lower bound for the left handside of~\eqref{e:agmon-prelim}.
Concerning the right handside of~\eqref{e:agmon-prelim}, we write for $E \in [\min_{\M} V  ,\max_{\M} V]$
$$
\sup_{\Omega_\alpha^-} |V -E - |\d_x\phi|^2 + \qe | \leq 4(\nor{V}{\infty} +1)+1  = : C_{V}.
$$
We fix $\eps_0 = \eps_0(\delta, \alpha)$ such that $\nor{\qe}{\infty} \leq \frac{1}{2}\alpha^2\delta$ for all $\eps \in (0,\eps_0)$.
Coming back to~\eqref{e:agmon-prelim}, we have obtained for $\delta \in (0,1)$ and $\eps\leq \eps_0$,
\begin{align*}
\int_0^L  \eps^2 |\d_x(e^{\phi/\eps}\psi)|^2 +  \frac{1}{2}\alpha^2\delta \int_{\Omega_\alpha^+} e^{2\phi/\eps}|\psi|^2 \leq C_{V}  \int_{\Omega_\alpha^-}e^{2\phi/\eps}|\psi|^2 .
\end{align*}
This implies 
\begin{align}
\label{e:agmon-sauf-rhs}
\int_0^L  \eps^2 |\d_x(e^{\phi/\eps}\psi)|^2 +  \frac{1}{2}\alpha^2\delta \int_0^L  e^{2\phi/\eps}|\psi|^2 \leq (C_{V} + 1)  \int_{\Omega_\alpha^-}e^{2\phi/\eps}|\psi|^2 .
\end{align}
To conclude the proof, we now estimate the right handside in~\eqref{e:agmon-sauf-rhs}. We write $\Omega_\alpha^- = \big(\Omega_\alpha^- \cap [0,\xzero) \big)\sqcup \big(\Omega_\alpha^- \cap [\xzero,L]\big)$ and split the integral accordingly, using that $V$ is injective on each part. We now only consider the second term, the first one being treated similarly.
 Uniform continuity of $V^{-1}$ on the compact $[\xzero,L]$ implies the existence of $\alpha=\alpha(\delta) \in (0, 1]$ such that 
\begin{align}
\label{e:unif-cont-delta-bis}
\left( E \in \R , \quad  x, y \in \{z; E \leq V(z) \leq E + \alpha^2 \} \cap [\xzero ,L] \right) \implies |x-y| \leq \delta .
\end{align} 
As a consequence, we have for $x \in \Omega_\alpha^- \cap [\xzero,L]$ (and $V(x)\geq E$, otherwise $\phi(x)=0$ and the same estimate is true)
$$
\phi (x)= (1-\delta) d_{A,E}(x) = (1-\delta) \int_{x_+(E)}^x \sqrt{(V(s) - E)_+} ds \leq (1-\delta) (x-x_+(E)) \alpha \leq (1-\delta)\delta \alpha  ,
$$
using~\eqref{e:unif-cont-delta-bis} (where $x_+(E)\in K_{E}$ is the solution in $[\xzero,L]$ of $V(x_+(E))=E$).
Coming back to~\eqref{e:agmon-sauf-rhs}, we now have 
\begin{align*}
\int_0^L  \eps^2 |\d_x(e^{\phi/\eps}\psi)|^2 +  \frac{1}{2}\alpha^2\delta \int_0^L  e^{2\phi/\eps}|\psi|^2 \leq (C_{V} + 1) e^{2\delta \alpha /\eps} \int_{\Omega_\alpha^-}|\psi|^2 \leq  (C_{V} + 1) e^{2\delta \alpha /\eps}  .
\end{align*}
We now want to replace $\phi$ by $d_{A,E}$. Recall that $\phi = (1-\delta) d_{A,E}$,  and that $0 \leq d_{A,E}(x) \leq L D_{V}$ for another constant $D_{V}:=  \sqrt{\max_{\M} V- \min_{\M} V + 1}$ uniformly in $x, E$, so that we may write 
\begin{align*}
& \int_0^L |\d_x(e^{d_{A,E}/\eps}\psi)|^2+  \int_0^L  e^{2d_{A,E}/\eps}|\psi|^2 \\
&\quad = \int_0^L |\d_x(e^{\delta d_{A,E}/\eps} e^{\phi/\eps} \psi)|^2+ \int_0^L |e^{\delta d_{A,E}/\eps} e^{\phi/\eps} \psi|^2\\
& \quad  \leq \left(1 + \frac{\delta D_{V}}{\eps}+\frac{2}{\alpha^2\delta}\right)e^{\delta L D_{V}/\eps}\left[ \int_0^L |\d_x(e^{\phi/\eps} \psi)|^2 + \frac{1}{2}\alpha^2\delta\int_0^L  e^{2\phi/\eps}|\psi|^2\right].
\end{align*}
Combining the above two estimates implies 
\begin{align*}
\int_0^L  \eps^2 |\d_x(e^{d_{A,E}/\eps}\psi)|^2 +  \int_0^L  e^{2d_{A,E}/\eps}|\psi|^2 \leq (C_{V} + 1)\left(1 + \frac{\delta D_{V}}{\eps}+\frac{2}{\alpha^2\delta}\right) e^{\frac{\delta}{\eps}\left( 2+ LD_{V}\right)}  ,
\end{align*}
which proves~\eqref{e:agmon-H1} up to changing $\delta \left( 2+ LD_{V}\right)$ into $\delta$. 

To obtain the bound on the normal trace, we need an $H^2$ bound on $e^{d_{A,E}/\eps}\psi$. To this aim, we follow e.g.~\cite[Remark~3.3]{Helffer:booksemiclassic} and first regularize $d_{A,E}$. We consider $\rho_\delta  = \frac{1}{\delta} \rho(\frac{\cdot}{\delta})\in C^\infty_c(-\delta,\delta)$ a nonnegative smooth approximation of the identity, and define $d_{A,E}^\delta = \rho_\delta \ast d_{A,E}$ for $\delta$ small enough, where $V, \qe$ (and $d_{A,E}$ accordingly) have been extended in a fixed neighborhood of $[0,L]$. We have $0\leq d_{A,E}^\delta  \leq \sup_{x\in [-\delta, L+\delta]}d_{A,E}(x) \leq 2LD_{V}$, and, uniformly for $x \in [0,L]$,  
\begin{align}
\label{e:daedelta}
|d_{A,E}^\delta (x) - d_{A,E}(x) |  & \leq \int |d_{A,E}(x- y) - d_{A,E}(x)| \rho_\delta(y) dy \leq D_{V} \int |y| \rho_\delta(y) dy \nonumber \\
&\leq  \delta D_{V} \int |y| \rho(y) dy , 
\end{align}
where we used that $|d_{A,E}'| = \sqrt{(V-E)_+} \leq D_{V}$. As a consequence, from~\eqref{e:agmon-H1}, we now obtain for a constant $D_V$ depending only on $V$, for $\eps< \eps_0 = \eps_0(\delta)$, denoting $\nor{ f}{H^1_\eps}=\e\nor{f'}{L^{2}}+\nor{f}{L^{2}}$ and $\Psi_\eps = e^{d_{A,E}^\delta/\eps}\psi$
\begin{align}
\label{e:agmon-smoothened}
\nor{ \Psi_\eps}{H^1_\eps} \leq 2\nor{e^{\frac{1}{\eps}(d_{A,E}^\delta-d_{A,E})}}{W^{1,\infty}} \nor{ e^{d_{A,E}/\eps}\psi}{H^1_\eps} \leq C_\delta \eps^{-1}e^{D_{V} \frac{\delta}{\eps}}   e^{\frac{\delta}{\eps}}.
\end{align}
The function $\Psi_\eps$ is solution of  
$$ 
( P_\eps -E ) \Psi_\eps = - 2\eps (e^{d_{A,E}^\delta/\eps})' \eps \psi'  - \eps^2 (e^{d_{A,E}^\delta/\eps})'' \psi , \quad  \Psi_\eps(0) = \Psi_\eps(L) = 0 .
$$
According to the above inequality~\eqref{e:agmon-smoothened} and bounds on $d_{A,E}^\delta$, we obtain $\nor{\Psi_\eps''}{L^2} \leq C_\delta e^{(D_V+2) \frac{\delta}{\eps}}$ uniformly for $E \in [\min_{\M} V-1 ,\max_{\M} V +1]$ and $\e \leq \e_{0}(\delta)$. This together with~\eqref{e:agmon-smoothened} directly implies $|e^{d_{A,E}^\delta(0)/\eps}\psi'(0)| = |\Psi_\eps'(0)| \leq C_\delta e^{D_V \frac{\delta}{\eps}}$ (and similarly at $L$). 
Using again~\eqref{e:daedelta} finally replaces $e^{d_{A,E}^\delta(0)/\eps}$ by $e^{d_{A,E}(0)/\eps}$ in this estimate with an additional $e^{CD_V \frac{\delta}{\eps}}$ loss, and thus implies~\eqref{e:Agmon-boundary} (after having changed $C D_V \delta$ into $\delta$). This concludes the proof of the theorem.
\enp

\subsection{Lower estimates in the classically allowed region}
\label{s:classically-allowed-section}
In this section, we first deduce the following ``geometric control estimate'' from the description of semiclassical measures in Theorem~\ref{thmmeasurex}. We then give a proof of Theorem~\ref{thmmeasurex}, relying on technical statements for semiclassical measures for one-dimensional boundary-value problems, proved in Section~\ref{s:app-msc} below.

\begin{proposition}[Geometric control in the classically allowed region]
\label{l:geom-control}
Let $V \in C^1([0,L])$ satisfy Assumption~\ref{assumptions} and $\qe \to 0$ in $C^1([0,L])$.
Then for any family $(\lambda_\eps)_{\eps\in (0,1)}$, $\lambda_{\eps}\in \R$, converging to zero as $\eps\to 0^+$, for any $\nu >0$, there are constant $C, \eps_0>0$ such that for all $\y \in \M$, all $\eps \in (0,\eps_0]$, all $E \in \R$ and $\psi$ satisfying~\eqref{e:eignfct-E-bis}, 
 we have  
\begin{align}
\label{e:obs-eign}
& \nor{\psi}{L^2(U)} \geq C , \quad \text{with } U= (\y-\nu,\y+\nu) \cap [0,L] ,\quad  \text{ if } E \geq  V(\y) - \lambda_{\eps}   , \\
&\frac{\eps}{\sqrt{|E|+1}} |\psi' (0)| \geq C ,\quad  \text{ if } E \geq  V(0)+ \nu , \label{e:obs-eign-0}  \\
&\frac{\eps}{\sqrt{|E|+1}} |\psi' (L)| \geq C ,\quad  \text{ if } E \geq  V(L)+ \nu  . \label{e:obs-eign-L}
\end{align}
\end{proposition}

Some remarks are in order:
\begin{itemize}
\item Note that the lemma states ``observability  inequalities'' for eigenfunctions~\eqref{e:eignfct-E-bis} from a neighborhood of a point $\y$, assuming a ``geometric control condition'', which is here formulated as $E \geq V(\y)$ (internal case) or $E \geq  V(0)+ \nu$ (observation from the boundary $0$) or $E \geq  V(L)+ \nu$ (observation from the boundary $L$). The latter condition ensures that all classical trajectories with energy $E$ intersect the region $(\y-\nu , \y+\nu)$ (internal case) or $0$ (observation from the boundary $0$) or $L$ (observation from the boundary $L$).
\item 
Note that the proof below proceeds by contradiction and uses semiclassical measures, following the general strategy introduced by Lebeau~\cite{Leb:96}. 
\item 
Note that the explicit expression of the measures in Theorem~\ref{thmmeasurex} can be used to describe for instance the asymptotic values of the constants $C$ in~\eqref{e:obs-eign}-\eqref{e:obs-eign-0}-\eqref{e:obs-eign-L}. Note also that the eigenfunction equation~\eqref{e:eignfct-E-bis} can be ``relaxed'' to a quasimode equation as in the statement of Theorem~\ref{thmmeasurex}.

\end{itemize}

\bnp[Proof of Proposition~\ref{l:geom-control} from Theorem~\ref{thmmeasurex}]
We proceed to the proof by contradiction, following the strategy introduced by Lebeau~\cite{Leb:96}. Given $\lambda_\eps\to 0$ and $\nu >0$, if the statement of the lemma is not satisfied, the following holds: for all $n \in \N$, there exist $\y_n\in \M$, $\eps_n \in (0,\frac{1}{n}]$, $E_n \in \R$, $\psi_n$ satisfying~\eqref{e:eignfct-E-n} with $r_{n}=0$, together with 
\begin{align}
\label{e:eignfct-E-n-obs}
\quad \nor{\psi_n}{L^2(\y_n-\nu,\y_n+\nu)} < \frac{1}{n} \quad \text{ in case }E_n \geq V(y_{n})-\lambda_{\eps_n}, \\
\label{e:eignfct-E-n-obs-bord}
 \text{resp. }\frac{\eps_{n}}{\sqrt{|E_{n}|+1}}|\psi_n' (0)|< \frac{1}{n} \text{ in case }E_n \geq V(0)+\nu, \\
  \text{resp. }\frac{\eps_{n}}{\sqrt{|E_{n}|+1}}|\psi_n' (L)|< \frac{1}{n}  \text{ in case }E_n \geq V(L)+\nu . \nonumber
\end{align}
We may now extract from the sequence $(\y_n ,\eps_n, E_n , \psi_n)_{n \in \N}$ a subsequence (which we do not relabel, with a slight abuse of notation) such that 
\begin{align*}
\eps_n \to 0  , \quad \y_n \to \y_* \in \M ,  \\
E_n \to E_* \in [V(y_*) , + \infty] , \\
|\psi_n(x)|^2 dx \rightharpoonup \mathfrak{m}_{E_*} ,
\end{align*}
where the last convergence holds in the sense of weak$-*$ convergence of measure.
The measure $\mathfrak{m}_{E_*}$ is described explicitly in Theorem \ref{thmmeasurex}. Note that the assumptions yield $E_{*}\geq V(y_*)  \geq V(\xzero)=\min V$. This implies $y_{*}\in [x_{-}(E_{*}),x_{+}(E_{*})]$ and in particular, $\mathfrak{m}_{E_*}( (\y_*-\nu/2 , \y_*+\nu/2))>0 $ in all three cases of the definitions of $\mathfrak{m}_{E_*}$ in Theorem \ref{thmmeasurex}.

Remark also that dominated convergence in~\eqref{e:eignfct-E-n-obs} implies that 
\begin{align}
\label{e:eignfct-E-n-obs-bis}
\nor{\psi_n}{L^2(\y_*-2\nu/3 , \y_*+ 2\nu/3 )}  \to 0 , \quad \text{as } n \to + \infty .
\end{align}
We obtain a contradiction with $\mathfrak{m}_{E_*}( (\y_*-\nu/2 , \y_*+\nu/2 ))>0$ by taking a bump function $\varphi\in C^{0}_{c}((\y_*-2\nu/3 , \y_*+2\nu/3),[0,1])$ equals to one on $(\y_*-\nu/2 , \y_*+\nu/2)$ which yields 
\bna
\nor{\psi_n}{L^2(\y_*-2\nu/3 , \y_*+2\nu/3)}^{2} \geq \nor{\varphi\psi_n}{L^2(0,L)}^{2}\underset{n\to +\infty}{\longrightarrow} \int_{[0,L]} \varphi(x)^{2} d\mathfrak{m}_{E_*}(x)\geq \mathfrak{m}_{E_*}( (\y_*-\nu/2 , \y_*+\nu/2 ))>0,
\ena
and contradicts~\eqref{e:eignfct-E-n-obs-bis}.
This proves the internal observability estimate~\eqref{e:obs-eign} and we are now left to prove the boundary observability. We only treat the case at the left boundary $x=0$, that is to prove that \eqref{e:eignfct-E-n-obs-bord} gives a contradiction.

To this aim, we now consider the cases $E_* = + \infty$ and $E_* < + \infty$ separately. 
If $E_* < +\infty$, then Theorem \ref{thmmeasurex} gives $|\eps_n\psi_n'(0)|^2 \to  2 C_{E_*}\sqrt{E_*-V(0)} \mathds{1}_{V(0) <E_*} $. Moreover, taking limit in the second part of \eqref{e:eignfct-E-n-obs-bord} gives $E_* \geq V(0)+\nu$. This implies $2 C_{E_*}\sqrt{E_*-V(0)} \mathds{1}_{V(0) <E_*} >0$ and therefore $\lim_{n\to +\infty}|\eps_n\psi_n'(0)|^2>0$ which is a contradiction to \eqref{e:eignfct-E-n-obs-bord}.

If now $E_* = +\infty$, Theorem \ref{thmmeasurex} gives $E_{n}^{-1}|\eps_n\psi_n'(0)|^2 \to  \frac{2}{L}$. Yet, since $E_{n}^{-1}\leq  \frac{2}{|E_{n}|+1}$ for $n$ large, \eqref{e:eignfct-E-n-obs-bord} gives $E_{n}^{-1}|\eps_n\psi_n'(0)|^2 \to  0$, which is a contradiction and ends the proof of the Lemma.
\enp

We are now left to prove Theorem \ref{thmmeasurex}. It relies on the following Proposition~\ref{propmeasureapp}  in which we describe fine localization properties and transport equations satisfied by semiclassical measures for solutions to $1D$ boundary value problems. The proof of Proposition~\ref{propmeasureapp} is given in Section~\ref{s:app-msc} below.
In the statement of Proposition~\ref{propmeasureapp}, we change slightly the current notation: we focus on the energy level $E=0$ for a potential $\V_{n}\to \V$, and consider the semiclassical parameter $h_n\to 0$. When deducing a proof of Theorem \ref{thmmeasurex}, we will use Proposition~\ref{propmeasureapp} both with 
\begin{itemize}
\item $h_n =\eps_n$ and $\V_{n}=V_n-E_n$ which converges to $\V=V-E_*$ (in case $E_n$ has a finite limit $E_*$), 
\item $h_n = \frac{\eps_n}{\sqrt{E_n}}$ and $\V_{n}=- 1+  \frac{V}{E_n}+\frac{\qen}{E_n}$ which converges to $\V=-1$ (in case $E_n \to +\infty$),
\end{itemize}
and in both cases, we describe the energy level $\V = 0$.
\begin{proposition}
\label{propmeasureapp}
Let $\V_{n}, \V\in C^1(\M)$ real valued so that $\V_{n}\to \V$ in $C^{1}(\M)$.  
Let  $h_{n}\to 0$ and $\psi_{n}$ be such that 
\begin{align}
\label{e:asspt-psi}
\psi_n \in H^2(0,L) \subset C^1([0,L]),  \quad \psi_{n}(0)=\psi_{n}(L)=0, \quad \nor{\psi_{n}}{L^{2}(0,L)}=1 , \nonumber \\
\quad -h_{n}^{2}\psi_{n}''+\V_{n}(x)\psi_{n}= r_n \quad \text{in }\mathcal{D}'((0,L)) ,
\end{align}
and, given a function $u$ defined on $[0,L]$, denote by $\udl{u}$ the function such that $\udl{u} = u$ on $[0,L]$ and $\udl{u}=0$ on $[0,L]^c$.
Assume that $r_n= \O_{L^2(0,L)}(h_n)$, then there exist 
\begin{itemize}
\item  a subsequence of indices (still denoted by $n$),
\item a probability measure $\mu$ on $T^*\R = \R_x\times\R_\xi$, supported in $[0,L]\times \R_{\xi}$, such that 
\begin{align}
\label{e:mu-exists}
\left( \Op_{h_n}(a) \udl{\psi_n} ,  \udl{\psi_n} \right)_{L^2} \to \langle \mu , a \rangle , \quad \text{ for all }a \in C^\infty_c(T^*\R),
\end{align}
\item two nonnegative numbers $\ell_{0}$ and $\ell_{L}$ so that
\bnan
\label{limtrathm}
|h_n\psi_n'(0^{+})|^{2}\to \ell_{0}, \quad |h_n\psi_n'(L^{-})|^{2}\to \ell_{L},
\enan
\item a probability measure $\mathfrak{m}$ on $\R$ such that $|\udl{\psi_n}(x)|^2 dx \rightharpoonup \mathfrak{m}$, in the sense of weak$-*$ convergence of measures on $\R$.
\end{itemize}
Moreover, writing $p(x,\xi) := \xi^2 + \V(x)$, the following statements hold:
\begin{enumerate}
\item  we have $\supp(\mu) \subset \{p(x,\xi)=0\}\cap [0,L]\times \R_{\xi}$;
\item if $r_n= \petito{h_n}_{L^2(0,L)}$, then $\mu$ satisfies $H_{p}\mu=0$ in $\mathcal{D}' \big((0,L)\times \R_{\xi} \big)$;
\item if $r_n= \petito{h_n}_{L^2(0,L)}$, then depending on the value $\V(0)$, we have  
 \begin{itemize}
 \item (Elliptic case): if $\V(0)>0$: then $\ell_{0}=0$ and there is $\delta >0$ such that $\mu=0$ in $(-\delta,\delta)\times \R$,
\item (Glancing case): if $\V(0)=0$: then $H_{p}\mu=-\ell_{0}\delta_{x=0}\otimes \delta'_{\xi=0} $ for $x$ close to $0$. If moreover \rouge{$\V'(0)\leq 0$,} then $\ell_{0}=0$ and $H_{p}\mu=0$ in $ \mathcal{D}' \big((-\infty,L)\times \R \big)$, 
\item (Hyperbolic case): if $\V(0)<0$: then $H_{p}\mu=\frac{\ell_{0}}{2\sqrt{-\V(0)}}\delta_{x=0}\otimes (\delta_{\xi=\sqrt{-\V(0)}}- \delta_{\xi=-\sqrt{-\V(0)}})$  in $ \mathcal{D}' \big((-\infty,L)\times \R \big)$, 
\end{itemize}
(and symmetric relations are true close to $L$).
\item The measures $\mathfrak{m}$ and $\mu$ are linked by $\mathfrak{m}=\pi^{*}\mu$, where $\pi : \R_x \times \R_\xi \to \R_x$ is the canonical projection, that is to say $\int_{\R}\varphi(x)~d\mathfrak{m}= \int_{\R^{2}}\varphi\circ \pi~d\mu$ for all $\varphi \in C^0_c(\R^2)$.
\end{enumerate}
\end{proposition}
 Note that since $H_{p}$ is only assumed to be a $C^{0}$ vector field, $H_{p}\mu$ is defined by duality which makes sense since $\mu$ is a measure (and not only a distribution), see Lemma~\ref{lmsupHpmu} below.
 
 Let us now prove Theorem \ref{thmmeasurex} from Proposition~\ref{propmeasureapp}. Note that the regularity assumption on $V \in C^1([0,L])$ requires some care in the propagation estimates for semiclassical measures (in the proof of Theorem \ref{thmmeasurex} as well as in the proof of Proposition~\ref{propmeasureapp}). One reason for this is that the Cauchy-Lipschitz theorem does not apply to the continuous Hamiltonian vector field $2\xi\d_x-V'(x)\d_\xi$. 
\bnp[Proof of Theorem \ref{thmmeasurex}]
We consider the cases $E_* = + \infty$ and $E_* < + \infty$ separately. In each case, we will compute a semiclassical measure, but with respect to a different small parameter namely $h_n = \frac{\eps_n}{\sqrt{E_n}}$ or $h_{n}=\e_{n}$ respectively. In the present proof, we shall describe the full semiclassical measure $\mu_{E_*}$ in phase space, associated to the sequence of eigenfunctions $\psi_n$ (extended by zero outside $[0,L]$) and the scale $h_n$. Then, the measure $\mathfrak{m}_{E_{*}}$ will be the (restriction to $[0,L]$ of the) projection in $x$ of the semiclassical measure $\pi^{*} \mu_{E_*} = \mathds{1}_{[0,L]} \mathfrak{m}_{E_{*}}$ , where $\pi : \R_x \times \R_\xi \to \R_x$ is the canonical projection (see the last Item of Proposition \ref{propmeasureapp}). The limits of the respective boundary terms will result from the computation of $\ell_{0}$ and $\ell_{L}$ in the same Proposition \ref{propmeasureapp}.
\paragraph{Case 1: $E_* = + \infty$.}
We rewrite the first equation in~\eqref{e:eignfct-E-n} as 
\begin{align*}
(-h_n^2\d_x^2 +\V_{n}  ) \psi_n =r_{n}E_{n}^{-1},
\end{align*}
where we have set $h_n = \frac{\eps_n}{\sqrt{E_n}}\to 0^+$ and $\V_{n}=- 1+  \frac{V}{E_n}+ \frac{\qen}{E_n}$. 
Extending $\psi_{n}$ by $0$ outside of $\M$ (without changing notation),  
Proposition~\ref{propmeasureapp} can be applied with $\V_{n}=- 1+  \frac{V}{E_n}+\frac{\qen}{E_n}$ and $\V=-1$ with $\V_n\to \V$ in $C^1([0,L])$ and $\frac{r_{n}}{E_{n}}=\petito{h_n}$ since $r_n=\petito{\eps_n}$. It provides a semiclassical measure $\mu$ such that, up to a subsequence,
$$
\left( \Op_{h_n}(a) \psi_n , \psi_n \right)_{L^2} \to \langle \mu , a \rangle , \quad \text{ for all }a \in C^\infty_c(\R\times\R).
$$
Moreover, according to Proposition~\ref{propmeasureapp}, the measure $\mu$ is supported by $[0,L]_x \times \{\pm 1\}_\xi$ and locally invariant by the flow of the vector field $2\xi \d_x$ in $(0,L)_x \times \{\pm 1\}_\xi$, we necessarily have 
\bna
\mu = \theta_1 \frac{\mathds{1}_{[0,L]}dx}{L} \otimes \delta_{\xi=1} + \theta_2 \frac{\mathds{1}_{[0,L]}dx}{L} \otimes \delta_{\xi=-1} + \theta_3 \delta_{(0,1)} +  \theta_4 \delta_{(0,-1)} + \theta_5 \delta_{(L,1)} +  \theta_6 \delta_{(L,-1)} , \\
\text{with }  \theta_j \in [0,1] , \quad \sum_j \theta_j=1.
\ena
(see below for a justification of this decomposition in a slightly more intricate setting).
Also, the second part of Proposition~\ref{propmeasureapp} gives
\begin{align}
\label{e:mumu}
2\xi \d_x \mu = \left(\frac{\ell_0}{2} \delta_{x=0}- \frac{\ell_L}{2} \delta_{x=L}\right) \otimes (\delta_{\xi=1}- \delta_{\xi=-1})  \quad \text{ on } \R_{x} \times \R_{\xi},
\end{align}
where $\ell_0$ and $\ell_L$ are the limits of the normal traces $h_n^2|\psi_n'(0)|^2$ and $h_n^2|\psi_n'(L)|^2$ respectively.

In particular, the derivative of $\mu$ is a measure. This implies that $\mu(\{(0,1)\}) = \mu(\{(0,-1)\})=\mu(\{(L,1)\})= \mu(\{(L,-1)\})=0$, and thus $\theta_j$ above vanish for all $j \geq 3$.
Therefore, there is $\theta \in [0,1]$ such that 
\bna
\mu = \theta \frac{\mathds{1}_{[0,L]}dx}{L} \otimes \delta_{\xi=1} + (1- \theta ) \frac{\mathds{1}_{[0,L]}dx}{L} \otimes \delta_{\xi=-1}   \quad \text{ on } \R \times \R .
\ena
Now, we compute the derivative of this measure, namely
$$
2\xi \d_x\mu = \frac{2\theta}{L} (\delta_{x=0} -\delta_{x=L})\otimes \delta_{\xi=1}  - \frac{2(1-\theta)}{L} (\delta_{x=0} -\delta_{x=L})\otimes \delta_{\xi=-1}.
$$ 
Identifying this with~\eqref{e:mumu} yields 
\begin{align}
\label{e:theta-ell-ell}
\theta = \frac{1}{2} ,\quad \text{and} \quad  \ell_0 = \ell_L = \frac{2}{L} .
\end{align}
We can now finally compute $\pi^{*} \mu=\frac{\mathds{1}_{[0,L]}dx}{L} $ which gives $\mathfrak{m}_{E_{*}}=\frac{dx}{L}$ after restriction to $[0,L]$. Since the limit is the same for any subsequence, we deduce that the convergence holds for the full sequence. Recalling that $h_n = \frac{\eps_n}{\sqrt{E_n}}$, the values of $\ell_0$ and $\ell_L$ in \eqref{e:theta-ell-ell} and the convergence result of \eqref{limtrathm} gives the expected limit for the boundary terms.
\paragraph{Case 2: $V(\xzero) \leq E_* < + \infty$.}

This time, we consider semiclassical operators scaled with the small parameter $h_{n}=\eps_n \to 0^+$, namely for $a \in C^\infty_c(\R_x\times \R_\xi)$, $\Op_{\eps_n}(a) = a(x, \eps_n D_x)$. 

Proposition \ref{propmeasureapp} applied with $h_{n}=\e_{n}$, $\V_{n}=V-E_{n}$ and $\V=V-E_{*}$ gives again a subsequence of indices (still denoted by $n$) and a nonnegative Radon measure $\mu$ on $T^*\R = \R_x \times \R_\xi$ such that 
$$
\left( \Op_{\eps_n}(a) \psi_n , \psi_n \right)_{L^2} \to \langle \mu , a \rangle , \quad \text{ for all }a \in C^\infty_c(\R\times \R) .  
$$
where we have again extended $\psi_{n}$ by zero without changing names.

Writing $p(x,\xi) = \xi^2 + V(x)$, Proposition~\ref{propmeasureapp} gives that $\mu$ is a probability measure, supported by the compact set 
\bna
p^{-1}(E_*) = \{ (x, \xi) \in [0,L] \times\R  \text{ such that } p(x,\xi) = E_* \} ,
\ena
and moreover invariant by the flow of the Hamiltonian vector field of $p$, namely $H_p = 2\xi \d_x - V'(x)\d_\xi$, locally in the interior of $(0,L)_x \times \R_\xi$. 
Note that, as already mentioned, we have slightly changed by a constant the notation for $p$ with respect to Proposition \ref{propmeasureapp} without changing the Hamiltonian flow. 

\medskip
We assume further in the proof that  
\begin{align}
\label{item:A3}
V(L) < V(0).
\end{align} 
The case $V(L) > V(0)$ is treated similarly. In the case $V(L) = V(0)$, there are actually less sub-cases to consider and the additional sub-case $E_* =V(L) = V(0)$ is treated as in sub-case 2 below (glancing near both endpoints of the interval); the two closed trajectories at energy $E_*$ are smooth and tangent to both boundaries $x=0$ and $x=L$.
 Given this additional assumption~\eqref{item:A3} on the shape of $V$, we only have to consider separately the following six sub-cases: 
\begin{enumerate}
\item $V(0) <E_* <+ \infty$
\item $E_* = V(0)$ 
\item $V(L) < E_* < V(0)$
\item $E_* = V(L)$ 
\item $V(\xzero) < E_* < V(L)$
\item $E_* = V(\xzero)$.
\end{enumerate}
\noindent
\paragraph{Sub-case 1: $V(0) <E_* <+ \infty$.}
Both $0$ and $L$ belong to $K_{E_*} = \pi(p^{-1}(E_*))$ (where $\pi : \R_x\times \R_\xi \to \R_x$ is the canonical projection)
 and the set $p^{-1}(E_*)$ decomposes as $p^{-1}(E_*) = \calC_+ \sqcup \calC_- \sqcup \{(0,\pm \sqrt{E_*-V(0)})\}\sqcup \{(L,\pm\sqrt{E_*-V(L)})\}$ where $\calC_\pm = \left\{(x, \pm \sqrt{E_*-V(x)}), x \in (0,L) \right\}$ are two disjoint bounded curves (that are both orbits of $H_p$ in case $V$ is regular enough). 
We may decompose accordingly the measure $\mu$ as 
\begin{align}
\label{e:decomp-mes-mes}
\mu = \mu \mathds{1}_{\calC_+} + \mu \mathds{1}_{\calC_-} + \mu \mathds{1}_{\{(0,\sqrt{E_*-V(0)})\}} +  \mu \mathds{1}_{\{(0,-\sqrt{E_*-V(0)})\}}+ \mu \mathds{1}_{\{(L,\sqrt{E_*-V(L)})\}} + \mu \mathds{1}_{\{(L,-\sqrt{E_*-V(L)})\}} ,
\end{align}
 in the sense of measures, i.e. for $F,E$ two Borel sets, $\mu \mathds{1}_{E}  (F) = \mu(E\cap F)$. In this decomposition, the four measures supported by points are proportional to Dirac masses.
We define $\delta_{\calC_\pm}$ as 
\begin{align}
\label{e:measuresCpm}
\left< \delta_{\calC_\pm} ,\varphi\right>=C_{E_*} \int_{0}^{L}\varphi(x,  \pm\sqrt{E_*-V(x)})\frac{dx}{\sqrt{E_*-V(x)}} , \quad  \text{with }  C_{E_*} = \left( \int_0^L(E_*-V(s))^{-\frac12} ds\right)^{-1}  ,
\end{align}
 for $\varphi\in C^0_c(\R_x\times \R_\xi)$ or, with a somewhat loose notation $$\delta_{\calC_\pm} = C_{E_*} \frac{\mathds{1}_{(0,L)}(x)dx}{\sqrt{E_*-V(x)}}\otimes \delta_{\xi = \pm\sqrt{E_*-V(x)}} .$$
 Let us now prove, using invariance by $H_p$, that $\mu \mathds{1}_{\calC_\pm}$ is proportional to $\delta_{\calC_\pm}$, that is, it is the unique invariant measure on $\calC_\pm$. This would be straightforward if we would have $V' \in C^1$, as a consequence of the Cauchy-Lipschitz theorem, but we only assume $V'\in C^0$ here.  
We define a measure $\nu$ on $(0,L)$ by
\begin{align*}
\left<\nu,f \right>_{\mathcal{M},C^0_c(0,L)}:=\left<\mu \mathds{1}_{\calC_+},\sqrt{E_*-V(x)}f\otimes 1 \right>_{\mathcal{M},C^0_c((0,L)\times \R)}, \quad \text{with } f\otimes 1 (x,\xi)=f(x)
\end{align*}
Let us first prove that $\d_x \nu=0$ in the distributional sense: we have
\begin{align*}
\left<\nu,\d_x f \right>_{\mathcal{M},C^0_c(0,L)}&=\left<\mu\mathds{1}_{\calC_+},\sqrt{E_*-V}\d_x f\otimes 1 \right>_{\mathcal{M},C^0_c((0,L)\times \R)}\\
&=\left<\mu\mathds{1}_{\calC_+},\xi\d_x f\otimes 1 \right>_{\mathcal{M},C^0_c((0,L)\times \R)}\\
&=\left<\mu\mathds{1}_{\calC_+},H_p (f\otimes 1) \right>_{\mathcal{M},C^0_c((0,L)\times \R)}=0
\end{align*}
where we have used that $\xi=\sqrt{E_*-V(x)}$ on $\calC_+$ in the first equality and the invariance property of $\mu\mathds{1}_{\calC_+}$ in the last equality (the latter is a consequence of the invariance of $\mu$ and the fact that $\overline{\calC_+} \cap \overline{\calC_-}=\emptyset$ in this subcase). This proves in particular that there exists a constant $\beta$ so that $\nu=\beta dx$. In particular, for $\varphi \in C^0_c((0,L)\times \R)$, we compute, using again that $\xi=\sqrt{E_*-V(x)}$ on $\calC_+$
\begin{align*}
\left<\mu\mathds{1}_{\calC_+},\varphi(x,\xi) \right>_{\mathcal{M},C^0_c((0,L)\times \R)}&=\left<\mu\mathds{1}_{\calC_+},\varphi(x,\sqrt{E_*-V(x)})\otimes 1 \right>_{\mathcal{M},C^0_c((0,L)\times \R)}\\
&=\left<\nu,\frac{\varphi(x,\sqrt{E_*-V(x)})}{\sqrt{E_*-V(x)}} \right>_{\mathcal{M},C^0_c(0,L)}=\beta\int_0^L \frac{\varphi(x,\sqrt{E_*-V(x)})}{\sqrt{E_*-V(x)}}~dx\\
&=\beta C_{E_*}^{-1}\left< \delta_{\calC_+} ,\varphi\right>.
\end{align*}

Coming back to the decomposition~\eqref{e:decomp-mes-mes} we have now obtained 
\bna
\mu = \theta_1 \delta_{\calC_+} + \theta_2  \delta_{\calC_-}+ \theta_3 \delta_{(0,\sqrt{E_*-V(0)})} +  \theta_4 \delta_{(0,-\sqrt{E_*-V(0)})} + \theta_5 \delta_{(L,\sqrt{E_*-V(L)})} +  \theta_6 \delta_{(L,-\sqrt{E_*-V(L)})} , \\
 \theta_j \in [0,1] , \quad \sum_j \theta_j=1.
\ena
Note also that for any $\varphi \in C^1_c(\R^2)$,
\begin{align}
\left< \delta_{\calC_\pm} ,H_{p}\varphi\right> &= C_{E_*} \int_{0}^{L}\left( \pm 2\sqrt{E_*-V(x)} \d_{x}-V'(x)\d_{\xi}\right)\varphi(x,  \pm\sqrt{E_*-V(x)})\frac{dx}{\sqrt{E_*-V(x)}} \nonumber \\
&= \pm 2 C_{E_*} \int_{0}^{L}\frac{d}{d x}\left[\varphi(x,  \pm\sqrt{E_*-V(x)})\right]dx \nonumber \\
& =\pm 2C_{E_*}\varphi(L,  \pm\sqrt{E_*-V(L)})\mp 2C_{E_*}\varphi(0,  \pm\sqrt{E_*-V(0)}).
\label{e:invar-Hpmu}
\end{align}
So that 
\bnan
\label{derivdC}
H_{p}\delta_{\calC_\pm} = \pm 2C_{E_*}\delta_{(0,\pm \sqrt{E_*-V(0)})} \mp 2C_{E_*}  \delta_{(L,\pm \sqrt{E_*-V(L)})} .
\enan
Moreover, both boundary points are of hyperbolic type (as in the case $E_* =\infty$). Using Proposition~\ref{propmeasureapp}, these measures satisfy the following equation 
\begin{align}
\label{eqnHpmubord}
H_{p} \mu & =\frac{\ell_0}{2\sqrt{E_*-V(0)}} \delta_{x=0}\otimes (\delta_{\xi=\sqrt{E_*-V(0)}}- \delta_{\xi=-\sqrt{E_*-V(0)}}) \nonumber\\
& \quad - \frac{\ell_L}{2\sqrt{E_*-V(L)}} \delta_{x=L}\otimes (\delta_{\xi=\sqrt{E_*-V(L)}}- \delta_{\xi=-\sqrt{E_*-V(L)}}),
\end{align}
where $H_p\mu$ is well-defined as a distribution according to Lemma~\ref{lmsupHpmu} below (using that the coefficients of $H_p$ are continuous and $\mu$ is a measure).
We can thus conclude as in the case $E_* =\infty$ that $H_{p}\mu$ is a measure and therefore  $\theta_{3}=\theta_{4}=\theta_{5}=\theta_{6}=0$. Again, comparing \eqref{derivdC}, \eqref{eqnHpmubord} and $\mu = \theta_1 \delta_{\calC_+} + \theta_2  \delta_{\calC_-}$, we obtain $\theta_1=\theta_2 = \frac{\ell_0}{4 C_{E_*}\sqrt{E_*-V(0)}}= \frac{\ell_L}{4 C_{E_*}\sqrt{E_*-V(L)}}$. The fact that $\mu$ is a probability measure gives $\theta_{1}=\theta_{2}=\frac{1}{2}$. In particular, $\mu=\frac{1}{2}(\delta_{\calC_+}+\delta_{\calC_-})$, $\ell_{0}=2C_{E_*}\sqrt{E_*-V(0)}$ and $\ell_{L}=2C_{E_*}\sqrt{E_*-V(L)}$ which gives the expected result for $\mathfrak{m}_{E_*}=\pi_{*}\mu$ and the limits of the boundary derivatives.
\paragraph{Sub-case 2: $E_* = V(0)$.}
The set $p^{-1}(E_*)$ writes $p^{-1}(E_*) = \calC_+ \sqcup \calC_- \sqcup  \{(0,0),(L,\pm\sqrt{E_*-V(L)})\}$ where again $\calC_\pm = \left\{(x, \pm \sqrt{E_*-V(x)}), x \in (0,L) \right\}$ are two disjoint bounded curves (that are both orbits of $H_p$ in case $V$ is regular enough). 
As in the first sub-case, we have accordingly 
\bna
\mu = \theta_1 \delta_{\calC_+} + \theta_2  \delta_{\calC_-}+ \theta_3 \delta_{(0,0)}
+ \theta_4 \delta_{(L,\sqrt{E_*-V(L)})} +  \theta_5 \delta_{(L,-\sqrt{E_*-V(L)})} , \\
 \theta_j \in [0,1] , \quad \sum_j \theta_j=1,
\ena
where $\delta_{\calC_\pm}$ is the unique invariant measure carried by $\calC_\pm$ and given by~\eqref{e:measuresCpm}.
Note that in the present situation, the right boundary $x=L$ is of hyperbolic type whereas the left boundary point $x=0$ is of diffractive type.
Now, the second part of Proposition \ref{propmeasureapp} yields in this case the equation
\begin{align}
\label{e:measure-glancing}
(2\xi \d_x- V' \d_{\xi}) \mu & = 
- \frac{\ell_L}{2\sqrt{E_*-V(L)}} \delta_{x=L}\otimes (\delta_{\xi=\sqrt{E_*-V(L)}}- \delta_{\xi=-\sqrt{E_*-V(L)}}),
\end{align}
  in a neighborhood of $x=L$.
In particular, the derivative $(2\xi \d_x- V' \d_{\xi}) \mu$ is a measure near $x=L$. This implies as in the above cases that $\theta_4=\theta_5=0$ and thus 
\bnan
\label{f:mutheta2}
\mu = \theta_1 \delta_{\calC_+} + \theta_2  \delta_{\calC_-}+ \theta_3 \delta_{(0,0)}.
\enan 
Near $x=L$, we are as in the previous Case 1 and differentiating this expression  (i.e. away from $x=0$) yields, using \eqref{derivdC},
$$(2\xi \d_x- V' \d_{\xi}) \mu =-  2 C_{E_*}  \delta_{x=L}\otimes \left( \theta_1  \delta_{\xi = \sqrt{E_*-V(L)}}- \theta_2  \delta_{\xi = -\sqrt{E_*-V(L)}} \right) .
$$
(See again Lemma~\ref{lmsupHpmu} below for the meaning of the left handside.)
Identifying the above two lines, we obtain again $\theta_1=\theta_2 = \frac{\ell_L}{4 C_{E_*}\sqrt{E_*-V(L)}}$.

We now consider the diffractive boundary at $x=0$ in~\eqref{e:measure-glancing}. Assumption~\ref{assumptions} and Proposition~\ref{propmeasureapp} imply that $H_{p}\mu=0$ close to $0$. A variant of \eqref{derivdC} implies $H_{p}\delta_{\calC_{\pm}}=0$ close to $0$, which combined with~\eqref{f:mutheta2} gives $\theta_3=0$. As a consequence $\theta_1= \theta_2 = \frac12$.

We have finally obtained that $\mu = \frac12 (\delta_{\calC_+} +  \delta_{\calC_-})$ and $\ell_L = 2 C_{E_*}\sqrt{E_*-V(L)}$. We can check that $\pi_{*}\mu$ gives the $\mathfrak{m}_{E_{*}}$ announced in the Theorem. It only remains to check that the Glancing case of Proposition~\ref {propmeasureapp} combined with the Assumption \ref{assumptions} (which implies $V'(0)<0$) impose $\ell_{0}=0$. This is the expected result since $ \sqrt{E_*-V(0)} \mathds{1}_{V(0) <E_*}=0$.  
\medskip
\noindent
\paragraph{Sub-case 3: $V(L) < E_* < V(0)$.} In this case, there is a single point $x_{E_*} \in (0,L)$ such that $V(x_{E_*})=E_*$ (it is given by $x_{E_*} =x_-(E_*)$), and we have $x_{E_*} <\xzero$  and $V'(x_{E_{*}})< 0$.
The set $p^{-1}(E_*)$ writes 
\begin{align} 
\label{e:decomposition-mes-C-delta}
p^{-1}(E_*) = \calC \sqcup  \{(L,\sqrt{E_*-V(L)})\}\sqcup  \{(L,- \sqrt{E_*-V(L)})\}, \\  \text{where} \quad \calC = \{(x, \pm \sqrt{E_*-V(x)}), x \in [x_{E_*},L)\}. 
\label{e:decomposition-mes-C-delta-bis}
\end{align}
We define the following probability measure on $\calC$ 
$$\left< \delta_{\calC} ,\varphi\right>= \frac{C_{E_*}}{2} \sum_{\pm} \int_{x_{E_*}}^{L}\varphi(x,  \pm\sqrt{E_*-V(x)})\frac{dx}{\sqrt{E_*-V(x)}},  \quad 
\text{with }C_{E_*}=\left(\int_{x_{E_*}}^L \frac{dx}{\sqrt{E_*-V(x)}} \right)^{-1} ,
$$
and we now aim at proving that $\mu \mathds{1}_{\calC}$ is proportional to $\delta_{\calC}$. Note that the difficulty in proving this comes again from the fact that $V'$ is only continuous. Would we have $V'\in W^{1,\infty}$, then the Cauchy-Lipschitz theorem would apply to $H_p$ and invariance of $\mu \mathds{1}_{\calC}$ would readily imply that it is proportional to $\delta_{\calC}$.
 
We define as above
\begin{align*}
\left< \delta_{\calC_\pm} ,\varphi\right>=C_{E_*} \int_{x_{E_*}}^L \varphi(x,  \pm\sqrt{E_*-V(x)})\frac{dx}{\sqrt{E_*-V(x)}} , 
\end{align*}
and we decompose further $\mu \mathds{1}_{\mathcal{C}} =\mu \mathds{1}_{\mathcal{C}_+}+\mu \mathds{1}_{\mathcal{C}_-} +\mu \mathds{1}_{\{(x_{E_*},0)\}}$. We notice that these measures are all compactly supported; we may test them with any function in $C^0(\R^2)$. The same proof as in Subcase 2 implies that necessarily  $\mu \mathds{1}_{\mathcal{C}_\pm} = \alpha_\pm \delta_{\mathcal{C}_\pm}$ and $\mu \mathds{1}_{\{(x_{E_*},0)\}}=\beta \delta_{(x_{E_*},0)}$.
Invariance of $\mu$ reads $\left< \mu ,H_{p}\phi\right> =0$ for all $\phi \in C^1(\R^2), \supp(\phi) \subset (0,L)\times \R$. Applied to $\phi(x,\xi)=\tilde\chi(x) \varphi(\xi)$ with $\tilde\chi \in C^1_c(0,L)$ such that $\tilde\chi(x_{E_*})=1$, we notice that $H_{p}\varphi =2 \xi \tilde\chi'(x)\varphi(\xi) -V'(x)\tilde\chi(x)\varphi'(\xi)$ and thus deduce 
$$
\langle \alpha_+ \delta_{\mathcal{C}_+}  + \alpha_- \delta_{\mathcal{C}_-} + \beta \delta_{(x_{E_*},0)} , 2 \xi \tilde\chi'(x)\varphi(\xi) -V'(x)\tilde\chi(x)\varphi'(\xi) \rangle = 0.
$$
Take $\chi \in C^\infty_c(\R)$ with $\chi =1$ in a neighborhood of zero and $\varphi_\epsilon(\xi) =\int_{-\infty}^{\xi/\epsilon} \chi(t) dt$. We obtain 
$$
0 = \langle \alpha_+ \delta_{\mathcal{C}_+}  + \alpha_- \delta_{\mathcal{C}_-} + \beta \delta_{(x_{E_*},0)} , 
2 \xi \tilde\chi'(x)\varphi_\eps(\xi) -V'(x)\tilde\chi(x)\frac{1}{\epsilon} \chi(\xi/\epsilon)
 \rangle ,
$$
whence multiplying by $\epsilon$, 
$$
0 = \epsilon \left\langle \alpha_+ \delta_{\mathcal{C}_+}  + \alpha_- \delta_{\mathcal{C}_-}  , 2 \xi \tilde\chi'(x)\varphi_\eps(\xi) -V'(x)\tilde\chi(x)\frac{1}{\epsilon} \chi(\xi/\epsilon)
 \right\rangle 
  -\beta V'(x_{E_*})   .
$$
Letting $\epsilon \to 0$ and using dominated convergence, we deduce 
$$
\beta V'(x_{E_*})  =  -\left\langle \alpha_+ \delta_{\mathcal{C}_+}  + \alpha_- \delta_{\mathcal{C}_-}  ,  V'(x)\tilde\chi(x)  \mathds{1}_{\{\xi=0\}}  \right\rangle =0
 $$
since $\mathcal{C}_\pm \cap \{\xi=0\}=\emptyset$. This implies $\beta V'(x_{E_*}) =0$, and thus $\beta=0$ since $V'(x_{E_*}) >0$.

Now we take any $\varphi \in C^1_c((0,L)\times \R)$ and compute as in~\eqref{e:invar-Hpmu}
$$
\left< \delta_{\calC_\pm} ,H_{p}\varphi\right> = C_{E_*}\int_{x_{E_*}}^{L} \pm\frac{d}{d x}\left[\varphi(x,  \pm\sqrt{E_*-V(x)})\right]dx=\mp C_{E_*} \varphi(x_{E_*},0).
$$ 
As a consequence, we obtain for all $\varphi \in C^1_c((0,L)\times \R)$
$$
0 = \left< \mu ,H_{p}\varphi\right> =\sum_\pm \alpha_{\pm}\left< \delta_{\calC_\pm} ,H_{p}\varphi\right>  = -\alpha_+ C_{E_*}\varphi(x_{E_*},0) + \alpha_-  C_{E_*}\varphi(x_{E_*},0) ,
$$
and thus $\alpha_+=\alpha_-$. This concludes the proof that $\mu \mathds{1}_{\calC}$ is proportional to $\delta_{\calC}$.

We now come back to the decomposition~\eqref{e:decomposition-mes-C-delta} and have obtained that 
\bna
\mu = \theta_1 \delta_{\calC} + \theta_2 \delta_{(L,\sqrt{E_*-V(L)})} +  \theta_3 \delta_{(L,-\sqrt{E_*-V(L)})} ,  \quad 
 \theta_j \in [0,1] , \quad \sum_j \theta_j=1.
\ena
The same computation as before gives 
$$
\left< \delta_{\calC} ,H_{p}\varphi\right> =  C_{E_*}\sum_{\pm} \int_{x_{E_*}}^{L} \pm\frac{d}{d x}\left[\varphi(x,  \pm\sqrt{E_*-V(x)})\right]dx=\sum_{\pm} \pm   C_{E_*} \varphi(L,  \pm\sqrt{E_*-V(L)}).
$$
So that $H_{p} \delta_{\calC}  = \sum_{\pm} \mp C_{E_*}\delta_{(L,\pm \sqrt{E_*-V(L)})}$ .
We thus argue as in the previous cases that $\theta_2=\theta_3=0$, hence $\theta_1=1$. As a consequence, $\mu = \delta_{\calC}$. Concerning the boundary estimates at $0$, we are in the Elliptic case of Proposition \ref{propmeasureapp} which implies $\ell_{0}=0$. This is the expected result since $ \sqrt{E_*-V(0)} \mathds{1}_{V(0) <E_*}=0$.  At $L$, we are in the Hyperbolic case, and we conclude as in the other subcases.

\medskip
\noindent
\paragraph{Sub-case 4: $E_*=V(L)$.}
The set $p^{-1}(E_*)$ writes $p^{-1}(E_*) = \calC \sqcup  \{(L,0)\}$ where $\calC \subset (0,L) \times \R$ is defined as in~\eqref{e:decomposition-mes-C-delta-bis} (and is an orbit of $H_p$ in case $V$ is regular). We have accordingly $\mu = \theta \delta_{\calC} + (1-\theta) \delta_{(L,0)}$, with $\delta_{\calC}$ the unique invariant measure carried by $\calC$ (a proof of uniqueness of this measure under the sole regularity assumption $V' \in C^0$ follows as in the above two sub-case). Moreover, the second part of Proposition \ref{propmeasureapp} yields in this case the equation 
\begin{align*}
(2\xi \d_x- V' \d_{\xi}) \mu & =0.
\end{align*}
The point $x=L$ is of diffractive type and the same analysis as in Sub-case 2 yields $\mu (\{(L,0)\}) =0$, hence $\theta=1$. This proves $\mu = \delta_\calC$ and we can conclude as in all other above cases. The proof that $\ell_{0}=\ell_{L}=0$ is performed as before for the respective Elliptic and Glancing cases (using that $V'(L)>0$).
\medskip
\noindent
\paragraph{Sub-case 5: $V(\xzero) < E_* < V(L)$.} The set $p^{-1}(E_*)$ is a $C^1$ closed curve contained in $(0,L)\times \R$ and $dp|_{p^{-1}(E_*)}$ does not vanish.
 The measure $\mu$ is supported on this curve and invariant by the vector field $H_p$, being nondegenerate and tangent to $p^{-1}(E_*)$. Henceforth, $\mu$ is the unique probability measure carried by $p^{-1}(E_*)$ and invariant by $H_p$ (again, uniqueness of this measure for $V' \in C^0$ follows as in the above sub-cases) defined by 
$$\left< \mu ,\varphi\right>=\frac{C_{E_*}}{2}\sum_{\pm} \int_{x_{-}(E_{*})}^{x_{+}(E_{*})}\varphi(x,  \pm\sqrt{E_*-V(x)})\frac{dx}{\sqrt{E_*-V(x)}},  \quad 
\text{with }C_{E_*}=\left(\int_{x_{-}(E_{*})}^{x_{+}(E_{*})} \frac{dx}{\sqrt{E_*-V(x)}} \right)^{-1}.
$$
The projection on $x$ of $\mu$ gives the expected result. Moreover, we are in the elliptic case at both boundaries $x=0$ and $L$ so that the normal trace converges to zero.
\medskip
\noindent
\paragraph{Sub-case 6: $E_* = V(\xzero)$.}
Note that the assumption $V(\y_*) \leq E_*$ implies that $V(\y_*)=E_*$ and thus $\y_* = \xzero$. 
We have $p^{-1}(E_*)= \{(\y_* , 0)\}$ and the only probability measure carried by this set is $\mu = \delta_{(x,\xi) = (\y_*,0)}$. We compute $\pi_{*}\mu=\delta_{\xzero}$ and we are again in the Elliptic case at both point of the boundary.

This concludes the proof of the theorem.
\enp

\subsection{Lower estimates in the classically forbidden region and near the turning points}
\label{s:class-forbid}
Next, we define the following ``semiclassical energy densities'' of the eigenfunctions $\psi$. For $x\in \M$: 
\begin{align*}
\E(x)&:= \eps^2 |\psi'|^2(x)+ |\psi|^2(x) , \\
 \E^+(x)&:=\eps^2 |\psi'|^2(x)+(V(x) -E )|\psi|^2(x) .
\end{align*}
The following lemma is a variant of~\cite[Lemma 12]{Allibert:98}, see also \cite[Lemma~4.10]{LL:18vanish}, in which we keep track of the dependence with respect to the lower order terms. 
\begin{lemma}[Tunneling into the classically forbidden region]
\label{l:tunneling-allbert-1D}
For all $\alpha>0$, all $E ,\psi ,\eps$ solution to~\eqref{e:eignfct-E-bis}, and all points $x,y$ belonging to the same connected component of $\{V-E \geq \alpha^2\}$, we have
\bna
\E^+(x)\leq \exp \left(\frac{2}{\e}\left| \int_x^y \sqrt{V(s) - E} ds \right| + \frac{\nor{V'}{\infty}}{\alpha^2}L + \frac{\nor{\qe}{\infty}}{\alpha\eps}L\right)\E^+(y)  .
\ena
\end{lemma}
\bnp[Proof of Lemma~\ref{l:tunneling-allbert-1D}]
We differentiate the function $\E^+$, yielding
$$(\E^+)'=2\eps^2\Re( \ovl{\psi}' \psi'')+ V'|\psi|^2 + 2(V-E) \Re(\psi \ovl{\psi}') .$$
We recall from the definition in~\eqref{e:def-Peps} and Equation~\eqref{e:eignfct-E-bis}  that we have
$$
E \psi =P_{\eps} \psi  = - \eps^2 \psi''  + V \psi  +  \qe \psi .
$$
This implies that
\begin{align}
\nonumber (\E^+)'& =2  \left(V -E + \qe  \right) \Re(\psi \ovl{\psi}') + V'|\psi|^2 + 2(V-E)\Re(\psi \ovl{\psi}') \\
\label{derivNRJ-E}& = (4(V-E) + 2 \qe)  \Re(\psi \ovl{\psi}')+ V'|\psi|^2  .
\end{align}
We now estimate each of the terms in the right handside of~\eqref{derivNRJ-E} on the set $\{V-E \geq \alpha^2\}$. We first have the pointwise estimate
\begin{align*}
\left| 4(V-E)\Re(\psi \ovl{\psi}') \right| & = 4 \e^{-1} \sqrt{V-E}  \left( \eps |\psi'| \right) \left(\sqrt{V-E} |\psi| \right) \\
 & \leq 2 \e^{-1} \sqrt{V-E} \left(\e^2 |\psi'|^2 + (V-E) |\psi|^2 \right) = 2\e^{-1} \sqrt{V-E}  \E^+ .
\end{align*}
Second, we have the pointwise estimate
 \begin{align*}
\left|V' |\psi|^2\right| = \frac{\left|V'\right|}{V-E}  (V-E) |\psi|^2 \leq \frac{\nor{V'}{\infty}}{\alpha^2}  \E^+ , \quad \text{ on } \{V-E \geq \alpha^2\}.
\end{align*}
 Third, we have on $\{V-E \geq \alpha^2\}$
$$
 \left|2 \qe  \Re(\psi \ovl{\psi}')\right|\leq \frac{\nor{\qe}{\infty}}{\eps} \left(\frac{\eps^2}{\alpha}|\psi'|^2 + \alpha |\psi^2| \right)
  \leq  \frac{\nor{\qe}{\infty}}{\eps}  \left(\frac{\eps^2}{\alpha}|\psi'|^2 + \alpha \frac{V-E}{\alpha^2}|\psi^2| \right)  =\frac{\nor{\qe}{\infty}}{\alpha \eps}  \E^+ .
$$ 
Combining the last three estimates in~\eqref{derivNRJ-E} yields, for all $t \in \{V-E \geq \alpha^2\}$
$$
\left| (\E^+)' (t)\right|\leq  \left( \frac{2}{\eps} \sqrt{V(t)-E} +  \frac{\nor{V'}{\infty}}{\alpha^2}  + \frac{\nor{\qe}{\infty}}{\alpha\eps} \right) \E^+(t) .
$$
Two applications of the Gronwall Lemma imply that for all $z<x$ such that $[z,x] \subset \{V-E \geq \alpha^2\}$, we have 
$$
e^{-\mu(x,z)} \E^+(z) \leq \E^+(x) \leq e^{\mu(x,z)}  \E^+(z) ,
$$
for $\mu(x,z) = \frac{2}{\eps} \int_{z}^x \sqrt{V(t)-E} dt+  \left(\frac{\nor{V'}{\infty}}{\alpha^2}+ \frac{\nor{\qe}{\infty}}{\alpha \eps }\right)(x-z)$. This yields the sought result.
\enp
Note that the estimate involving $\nor{V'}{\infty}$ could be slightly refined using a sign assumption on $V'$.

The following Lemma is an analogue of \cite[Lemma 11]{Allibert:98}, see also \cite[Lemma~4.11]{LL:18vanish}, and gives a rough Gronwall type estimate for the energy $\E$, without precise constants. 
The interest of this less precise result is that it remains true uniformly for all $x\in \M$. This allows in particular to compensate the fact that Lemma~\ref{l:tunneling-allbert-1D} is not uniform when $x$ is close to the boundary of the set $\{V-E>0\}$.

\begin{lemma}[Rough Gronwall estimate]
\label{l:transitions}
For all $E\in \R, \psi \in H^2(\M)\cap H^1_0(\M)$ and all $\eps>0$ such that $P_\eps \psi = E \psi$, and all $x,y \in \M$, we have
\bna
\E(x)\leq \exp \left(\frac{1}{\e}  |x-y|\left(  \nor{V - E+1}{L^\infty(I_{x,y})} + \nor{\qe}{\infty} \right) \right)\E(y)  ,
\ena
where $I_{x,y}$ is the interval between $x$ and $y$.
\end{lemma}

\bnp
The proof is very close to that of Lemma~\ref{l:tunneling-allbert-1D}. We write similarly
\begin{align*}
(\E)'& =2\eps^2\Re( \ovl{\psi}' \psi'')+ 2 \Re(\psi \ovl{\psi}')   =2(V-E + \qe + 1)\Re(\psi \ovl{\psi}') .
\end{align*}
This implies on the interval $I_{x,y}$
\begin{align*}
 |(\E)'|& \leq  \frac{1}{\eps}\left( \nor{V -E +1}{L^\infty(I_{x,y})} +  \nor{\qe}{\infty} \right) \E ,
\end{align*}
and we conclude the proof with a Gronwall argument on $I_{x,y}$ as in Lemma~\ref{l:tunneling-allbert-1D}.
\enp

\subsection{End of the proof of Theorem~\ref{t:allibert-lower-uniform}}
\label{l:end-of-proof}
With the three previous Lemmata at hand, we are now in position to prove Theorem~\ref{t:allibert-lower-uniform}. We first prove the following intermediate result.

\begin{lemma}[Lower bounds on eigenfunctions]
\label{l:allibert-lower-uniform}
Suppose $V,V_\eps$ satisfy Assumption~\ref{assumptions}.
Then, there is a constant $D>0$ such that for any $\mathbf{y}_0\in [0,L]$ and any $\delta >0$, there is $\eps_0>0$ such that for all 
$E \in \R$, $0<\eps<\eps_0$ and $\psi$ solution to~\eqref{e:eignfct-E-bis},
 we have, 
\begin{align}
&\nor{\psi}{L^2(U)} \geq e^{-\frac{1}{\eps}(d_{A,E}(\mathbf{y}_0) + D\delta)}  , \quad  U=(\mathbf{y}_0-\delta, \mathbf{y}_0+\delta) \cap [0,L] , 
\label{e:internal-obs-bis} \\
&\frac{\eps}{\sqrt{|E|+1}} |\psi'(0)|  \geq e^{-\frac{1}{\eps}(d_{A,E}(0) + \delta)}, \quad \frac{\eps}{\sqrt{|E|+1}} |\psi'(L)|  \geq e^{-\frac{1}{\eps}(d_{A,E}(L) + \delta)} \label{e:bound-obs-bis}  .
\end{align}
\end{lemma}
Note that in this statement (as well as in all statements of the article), $\delta$ is thought of as a small parameter.

\bnp[Proof that Lemma~\ref{l:allibert-lower-uniform} implies Theorem~\ref{t:allibert-lower-uniform}] 
Notice first that according to Remark~\ref{r:EgeqE0}, it suffices to consider $E\geq E_0$.
Then, the only difference between the two statements concerns the internal observation. We write $U=[\mathbf{z}_1,\mathbf{z}_2]$ with $\mathbf{z}_1,\mathbf{z}_2\in [0,L]$. We treat the case  for which $\mathbf{z}_1 \geq \xzero$: the case $\mathbf{z}_2 \leq \xzero$ is treated similarly. Concerning the case $\mathbf{z}_1 < \xzero <\mathbf{z}_2$, we take $\mathbf{y}_0=\xzero$ and choose $\delta >0$ small enough so that $(\xzero-\delta, \xzero+\delta)\subset (\mathbf{z}_1 ,\mathbf{z}_2)$ and Lemma~\ref{l:allibert-lower-uniform} yields the result since in this case $\inf_{x\in U}d_{A,E}(x)=d_{A,E}(\xzero)$ for all $E\geq E_0$.

Since we assume now $\mathbf{z}_1 \geq \xzero$, we have $\inf_{x\in U}d_{A,E}(x)=d_{A,E}(\mathbf{z}_1)$. 
According to Lemma~\ref{l:regularity-x-da}, $d_{A,E}$ is uniformly Lipschitz, so there $\widetilde{\delta}>0$ small enough uniform in $E\geq E_{0}$ so that $|d_{A,E}(z)-d_{A,E}(\mathbf{z}_1)|\leq \delta$ for $|z- \mathbf{z}_1|\leq \widetilde{\delta}$. We can also assume $\widetilde{\delta}<(\mathbf{z}_2-\mathbf{z}_1)/2$ and $\widetilde{\delta}\leq \delta$. Applying Lemma~\ref{l:allibert-lower-uniform} with $\mathbf{y}_0=\mathbf{z}_1+\widetilde{\delta}$ and $\delta$ replaced by $\widetilde{\delta}$, we obtain $\nor{\psi}{L^2((\mathbf{y}_0-\widetilde{\delta}, \mathbf{y}_0+\widetilde{\delta}) \cap [0,L] )} \geq e^{-\frac{1}{\eps}(d_{A,E}(\mathbf{y}_0) + D\widetilde{\delta})} $.  Since $(\mathbf{y}_0-\widetilde{\delta}, \mathbf{y}_0+\delta) \cap [0,L]\subset U$ and using the previous estimates, this gives  $\nor{\psi}{L^2(U)} \geq e^{-\frac{1}{\eps}(d_{A,E}(\mathbf{z}_1) + (D+1)\delta)} $ which is the expected result up to changing $\delta$.
\enp

We now prove Lemma~\ref{l:allibert-lower-uniform}, as a consequence of Theorem~\ref{thmmeasurex} and Lemmata~\ref{l:tunneling-allbert-1D} and~\ref{l:transitions}.
\bnp[Proof of Lemma~\ref{l:allibert-lower-uniform}]
We first prove the internal observation inequality~\eqref{e:internal-obs-bis}.
We distinguish different cases according to the respective location of the points $\mathbf{y}_0$ and $\xzero$.

\medskip
Consider first the case where $\xzero \in (\mathbf{y}_0-\delta , \mathbf{y}_0+\delta)$. Then, Proposition~\ref{l:geom-control} 
with $\nu$ small enough so that $(\mathbf{y}_0-\nu, \mathbf{y}_0+\nu) \subset U$,  yields 
$$
\nor{\psi}{L^2(U)} \geq \nor{\psi}{L^2(\xzero-\nu, \xzero+\nu)} \geq C_0 ,
$$
uniformly for $E \in \R$, which implies~\eqref{e:internal-obs-bis} in this case.

\medskip
We now consider the case where $\xzero \notin (\mathbf{y}_0-\delta , \mathbf{y}_0+\delta)$, and assume further in what follows that $\xzero \leq  \mathbf{y}_0 -\delta$. The case $\xzero \geq \mathbf{y}_0 +\delta$ is proved similarly (by symmetry). In particular, this implies $V'(\mathbf{y}_0) >0$ and $V(\mathbf{y}_0)>\min_{\M} V$.

 For this $\delta$, Proposition~\ref{l:geom-control} yields the existence of $C_0, \eps_0>0$ such that for all $z \in \M$, all $\eps \in (0,\eps_0]$, all
 $E \in \R$ and $\psi$ solution to~\eqref{e:eignfct-E-bis}, we have 
\begin{align}
\label{e:GCC-unif}
E \geq  V(z) \implies \nor{\psi}{L^2((z-\frac{\delta}{2},z+\frac{\delta}{2})\cap \M)} \geq C_0 . 
\end{align}
 Thanks to a variant of Remark \ref{r:EgeqE0}, we can assume from now on that $E\geq E_{0}$.

\paragraph{Case 1: $x_+(E)\geq \mathbf{y}_0 - \delta/2$.}
In this case, either $x_+(E)\geq \mathbf{y}_0$ (hence $E\geq V(\mathbf{y}_0)$), then~\eqref{e:GCC-unif} with $z=\mathbf{y}_0$ yields 
$$
\nor{\psi}{L^2\left((\mathbf{y}_0-\delta,\mathbf{y}_0+\delta)\cap \M\right)}\geq \nor{\psi}{L^2\left((\mathbf{y}_0-\delta/2,\mathbf{y}_0+\delta/2)\cap \M\right)} \geq  C_0 ,
$$
which concludes the proof in that case.

\medskip
Or else if $x_+(E)\leq \mathbf{y}_0 \leq x_+(E)+\delta/2$, then~\eqref{e:GCC-unif} with $z=x_+(E)$ yields 
$$
\nor{\psi}{L^2\left(((\mathbf{y}_0-\delta,\mathbf{y}_0+\delta)\cap [0,L]\right)}\geq \nor{\psi}{L^2\left(((x_+(E)-\frac{\delta}{2},x_+(E)+\frac{\delta}{2})\cap [0,L]\right)} \geq  C_0 ,
$$
which concludes the proof in that case.

\paragraph{Case 2: $x_+(E)< \mathbf{y}_0 - \delta/2$.} 
Lemma~\ref{l:regularity-x-da} (uniform continuity of $V^{-1}$ on the compact $[\xzero, L]$) implies the existence of $\alpha >0$ such that for all $x,y \in \M, E\in \R$,
\begin{align}
\label{e:unif-cont-delta}
    x, y \in \{z \in [\xzero,L] , E-\alpha^{2} \leq V(z) \leq E + \alpha^2 \} \implies |x-y| \leq \delta/4.
\end{align}

In this case, that $V(x_+(E))=E$ together with~\eqref{e:unif-cont-delta} implies that $\mathbf{y}_0 \notin  \{z \in [\xzero,L] , E-\alpha^{2} \leq V(z) \leq E + \alpha^2 \}$. Since $x_+(E)< \mathbf{y}_0$ in this case, this implies that necessarily $V( \mathbf{y}_0 )>E +\alpha^2$. 
  Estimate~\eqref{e:GCC-unif} with $z=x_+(E)$ implies
\bnan
\label{e:z0E}
C_0 \leq  \nor{\psi}{L^2\left((x_+(E)-\frac{\delta}{2},x_+(E)+\frac{\delta}{2})\cap \M\right)}  .
\enan
Lemma~\ref{l:transitions} together with $\nor{\qe}{\infty}\leq 1$ yields
\bnan
\label{e:trainsition-grongron}
|\psi|^2(x)\leq \exp \left(\frac{1}{\e}  |x-y| (2+2\|V\|_\infty) \right)\E(y)  , \quad x,y \in \M .
\enan
Integrating over $x \in (x_+(E)-\frac{\delta}{2} ,x_+(E)+\frac{\delta}{2})\cap [0,L]$ implies, for $y = x_+(E)+\frac{\delta}{2}< \mathbf{y}_0\leq L$, 
\bnan
\label{e:interm-allib}
 \nor{\psi}{L^2(x_+(E)-\frac{\delta}{2},x_+(E)+\frac{\delta}{2})\cap [0,L]}^2
 \leq \delta  \exp \left(\frac{\delta }{\e}(2+2\|V\|_\infty) \right)\E\left(x_+(E)+\frac{\delta}{2}\right) .
\enan
Now, remark that~\eqref{e:unif-cont-delta} implies $0 < x_+(E+ \alpha^2) -x_+(E) \leq \delta/4$.
The point $y = x_+(E)+\frac{\delta}{2} \in \{z;V(z)-E \geq \alpha^2\}$ is chosen so that to apply Lemma~\ref{l:tunneling-allbert-1D}.
Note first that on the set $\{z;V(z)-E \geq \alpha^2\}$ and for $|\alpha|<1$ (which we may assume), we have $\E \leq \alpha^{-2}\E^+$ and that $z \geq x_+(E)+\delta/4 \implies z \in \{V-E \geq \alpha^2\}$ (this is the case for $z= \mathbf{y}_0$).
 Lemma~\ref{l:tunneling-allbert-1D} now implies, for all $z \geq x_+(E)+\delta/4$,
\begin{align}
\label{e:allibert-ouvert-E}
\alpha^2 \E( x_+(E)+\delta/4 ) & \leq \E^+(x_+(E)+\delta/4)  \nonumber \\
& \leq \exp \left(\frac{2}{\e}\left| \int_{x_+(E)+\delta/4}^{z} \sqrt{V(s) - E} ds \right| + \frac{\nor{\qe}{\infty}}{\alpha\eps}L+ \frac{\nor{V'}{\infty}}{\alpha^2}L \right)\E^+(z)  .
\end{align}
Integrating in $z \in  (\mathbf{y}_0-\delta/4 ,\mathbf{y}_0 - \delta/8)$ (which implies $z \geq x_+(E)+\delta/4$ according to the assumption $x_+(E)< \mathbf{y}_0 - \delta/2$) yields
\bna
 \frac{\delta}{8} \alpha^2 \E( x_+(E)+\delta/2 ) \leq  \exp \left(\frac{2}{\e}\left| \int_{x_+(E)+ \delta/4}^{\mathbf{y}_0 } \sqrt{V(s) - E} ds \right| + \frac{\nor{\qe}{\infty}}{\alpha\eps}L+ \frac{\nor{V'}{\infty}}{\alpha^2}L \right)
\int_{\mathbf{y}_0-\delta/4}^{\mathbf{y}_0  - \delta/8} \E^+(s) ds.
\ena
An interpolation estimate together with $P_\eps \psi=E\psi$ yields
\begin{align*}
\int_{\mathbf{y}_0-\delta/4}^{\mathbf{y}_0  - \delta/8} \E^+(s) ds
&  \leq C\delta^{-1} \left( \nor{\psi}{L^2(\mathbf{y}_0-\delta/4,\mathbf{y}_0-\delta/8)}^2 + \nor{\psi}{L^2(\mathbf{y}_0-\delta/2 ,\mathbf{y}_0)}\nor{\eps^2 \psi''}{L^2(\mathbf{y}_0-\delta/2 ,\mathbf{y}_0)}\right) \\
&  \leq C\delta^{-1} \nor{\psi}{L^2(\mathbf{y}_0-\delta ,\mathbf{y}_0+\delta)\cap [0,L]}^2 .
\end{align*}
 Note that we have used $E\leq \nor{V}{\infty}$ otherwise this zone is empty.
Combining the above two estimates with~\eqref{e:z0E} and \eqref{e:interm-allib} yields the existence of constants $C =C(V, \delta,L)>0$ (recall that $\alpha$ depends on $\delta$ and $V$) independent on $E, \eps$ such that
\begin{align*}
1 &\leq C  \exp \left\{ \frac{2}{\e} \left( \left| \int_{x_+(E)+ \delta/4}^{\mathbf{y}_0} \sqrt{V(s) - E} ds \right| +(2+2\|V\|_\infty) \delta+ \frac{\nor{\qe}{\infty}}{\alpha}L \right)\right\}\nor{\psi}{L^2(\mathbf{y}_0-\delta ,\mathbf{y}_0+\delta)\cap [0,L]}^2.
\end{align*}
We further assume that $\eps_0$ is sufficiently small so that  assume that $\frac{\nor{\qe}{\infty}}{\alpha}L \leq \delta$ for all $\eps \in (0,\eps_0)$.
This then concludes the proof in that case, and hence the proof of~\eqref{e:internal-obs-bis} in the theorem.
 
\bigskip

We now explain how this proof needs to be modified in the case of boundary observability~\eqref{e:bound-obs-bis}, say, from the right boundary point $L$. In this case, the range of energy levels $E \in \R$ is again split in three different regimes.
We fix again $\alpha>0$ as in~\eqref{e:unif-cont-delta}.

First, if $E\geq V(L) +1$ then Proposition~\ref{l:geom-control} Estimate~\eqref{e:obs-eign-L} (taken for $\nu= 1$) yields $\frac{\eps}{\sqrt{|E|+1}}|\psi' (L)| \geq C$, which concludes the proof in that case.

Second, we consider the case $V(L)-\alpha^2 \leq E \leq V(L) +  1$. We remark that we have again, by definition of $\alpha$ and $x_+$, 
$$
 L-\delta/2 \leq x_+(V(L)-\alpha^2) \leq x_+(V(L)) = L .
$$
 Hence $(x_+(V(L)-\alpha^2)-\frac{\delta}{2},x_+(V(L)-\alpha^2)+\frac{\delta}{2})\cap \M \subset (L - \delta , L]$. Applying Estimate~\eqref{e:GCC-unif} for $z=x_+(V(L)-\alpha^2)$ and using $V(L)-\alpha^2 \leq E $, yields
$$
 \nor{\psi}{L^2(L-\delta , L)}  \geq C_0.
$$
Using \eqref{e:trainsition-grongron} integrated in $x \in (L-\delta , L)$ and taken for $y =L$ implies
$$
C_0^2 \leq  \nor{\psi}{L^2(L-\delta , L)}^2 \leq C \exp \left(\frac{\delta}{\e}  (2+2\|V\|_\infty)  \right)\E(L) ,
$$
where $\E(L) = \eps^2|\psi'(L)|^2$ on account to the Dirichlet boundary condition. This concludes the proof in that case.

Third, if $E\leq V(L)-\alpha^2$, the proof follows exactly as in the Case 2 above for the proof of~\eqref{e:internal-obs-bis}, except that the proof is finished when writing Estimate~\eqref{e:allibert-ouvert-E}, at the point $z=L$, together with noticing that $\E^+(L) = \eps^2|\psi'(L)|^2$, on account to the Dirichlet boundary condition.

This concludes the proof of~\eqref{e:bound-obs-bis} at the right boundary point $L$, and the proof is the same at the left boundary point $0$.
\enp

\section{Semiclassical measures for one-dimensional boundary-value problems}
\label{s:app-msc}
The object of this Section is to make precise different properties of semiclassical measures in the presence of boundary (and in dimension one only). 
 The combination of all results proved in the section constitutes a proof of Proposition~\ref{propmeasureapp}. 
The proof relies only on standard facts of semiclassical analysis for which we refer e.g. to~\cite{Robert:book,DS:book,Zworski:book} and semiclassical measures~\cite{Gerard:Bloch91,GL:93,GMMP:97,Zworski:book}.
Concerning the boundary value problem, we essentially follow~\cite{GL:93} with several major simplifications (due to absence of geometry of the boundary) and some minor complications (due to the family of limited regularity potentials converging in $C^1$).
We thus present a self-contained proof except for usual semiclassical analysis and semiclassical measure in $1D$. The latter material can be found in~\cite[Chapters~4 and~5]{Zworski:book} for instance.

To make the reading easier, we divide the proof in several Lemmata. 

\subsection{Regularity and traces}
\label{s:reg-of-traces}
We begin with standard regularity estimates (see e.g.~\cite[Lemma~2.1]{GL:93}). 
\begin{lemma}
\label{l:preliminary-trace} There is $C>0$ such that for all $h\in (0,1)$, $r \in L^2(0,L)$, $\mathscr{V} \in L^\infty(0,L)$ and $\psi \in H^2(0,L) \subset C^1([0,L])$ such that 
$$
\psi(0)=\psi(L)=0,  \quad -h^{2}\psi''+\mathscr{V} \psi= r \quad \text{in }\mathcal{D}'((0,L)) ,
$$
we have 
\begin{align}
\label{e:estim-bebete}
h^2 \|\psi'\|_{L^2(0,L)}^2 &  \leq \|\mathscr{V}\|_{L^\infty(0,L)}\|\psi\|_{L^2(0,L)}^2 + \|r\|_{L^2(0,L)} \|\psi\|_{L^2(0,L)}, 
\end{align}
and if moreover $\mathscr{V} =\mathscr{V}_{1} +\mathscr{V}_{2} $ with  $\mathscr{V}_{2} \in C^1([0,L])$ and $h \in(0,1)$,
\begin{align}
\label{e:estim-bebete-trace}
h^2|\psi'|^{2}(0^+)+ h^2|\psi'|^{2}(L^-) &\leq 
C\left(h^{-2}\|\mathscr{V}_{1}\|_{L^\infty(0,L)}^2+ \|\mathscr{V}_{2}\|_{C^1(0,L)} + 1\right) \|\psi\|_{L^2(0,L)}^2 
+Ch^{-2}\|r\|_{L^2(0,L)}^2  . 
\end{align}
\end{lemma}
Note that all along the present Section~\ref{s:app-msc}, we have $\mathscr{V}_{2}=\mathscr{V}_{2} \in C^1([0,L])$
\bnp
Multiplying the equation by $\ovl{\psi}$, integrating on $(0,L)$ and using an integration by parts, we obtain 
\bna
h^{2}\int_{(0,L)}|\psi'|^{2}dx+\int_{(0,L)}\mathscr{V}(x)|\psi|^{2}dx=\int_{(0,L)}r\ovl{\psi}dx.
\ena
The Cauchy-Schwarz inequality yields~\eqref{e:estim-bebete}.
To prove the second inequality, multiply the equation by $\chi(x)\ovl{\psi}'$ with $\chi \in C^\infty_c(\R ; [0,1])$ equal to $-1$ near $0$ and equal to $1$ near $L$. Integrating, we obtain 
\begin{align}
\label{e:toto-titi-tata}
0&=h^{2}\Re \int_{(0,L)}\psi''\chi \ovl{\psi}'dx - \Re \int_{(0,L)}\mathscr{V}(x)\psi\chi \ovl{\psi}'dx+\Re \int_{(0,L)}r\chi \ovl{\psi}'dx
\end{align}
Next integrating by parts, we obtain for the first term of~\eqref{e:toto-titi-tata}
$$
h^{2}\Re \int_{(0,L)}\psi''\chi \ovl{\psi}'dx =\frac{h^{2}}{2}\int_{(0,L)}\chi \frac{d}{dx}|\psi'|^{2}dx =\frac{h^{2}}{2}\left[|\psi'|^{2}(0^+)+|\psi'|^{2}(L^-)\right] -\frac{h^{2}}{2}\int_{(0,L)}\chi' |\psi'|^{2}dx.
$$
Concerning the last term of~\eqref{e:toto-titi-tata}, we simply write
$$
\left| \Re \int_{(0,L)}r\chi \ovl{\psi}'dx \right| \leq \|r\|_{L^2(0,L)} \|\psi'\|_{L^2(0,L)} \leq h^2 \|\psi'\|_{L^2(0,L)}^2 +  h^{-2}\|r\|_{L^2(0,L)}^{2} .
$$
We may estimate the second term of~\eqref{e:toto-titi-tata} with $\mathscr{V}=\mathscr{V}_1+\mathscr{V}_2$,  $\mathscr{V}_1\in L^\infty , \mathscr{V}_2 \in C^1$ as 
$$
\left| \Re \int_{(0,L)}\mathscr{V}_1(x)\psi\chi \ovl{\psi}'dx \right|  \leq \|\mathscr{V}_1\|_{L^\infty} \|\psi\|_{L^2} \|\psi'\|_{L^2} \leq  h^{-2}\|\mathscr{V}_1\|_{L^\infty}^2 \|\psi\|_{L^2}^2 +  h^2\|\psi'\|_{L^2}^2  ,
$$
and, integrating by parts, using $\psi(0) = \psi(L)=0$,  
\begin{align*}
\left| \Re \int_{(0,L)}\mathscr{V}_{2}(x)\psi\chi \ovl{\psi}'dx \right| =
\left|  \frac{1}{2} \int_{(0,L)}\mathscr{V}_{2}(x)\chi \frac{d}{dx}|\psi|^{2}dx \right| = \left| \frac{1}{2}\int_{(0,L)}(\mathscr{V}_{2}\chi)'|\psi|^{2}dx\right| \leq C\|\mathscr{V}_{2}\|_{C^1(0,L)}\|\psi\|_{L^2(0,L)}^2 .
\end{align*}
Combining the above four lines in~\eqref{e:toto-titi-tata} implies 
\begin{align*}
h^2|\psi'|^{2}(0^+)+ h^2|\psi'|^{2}(L^-) 
&  \leq Ch^2 \|\psi'\|_{L^2(0,L)}^2 + +C h^{-2}\|\mathscr{V}_{1}\|_{L^\infty}^2 \|\psi\|_{L^2}^2 \\
 & \quad +C\|\mathscr{V}_{2}\|_{C^1(0,L)}\|\psi\|_{L^2(0,L)}^2  +C h^{-2}\|r\|_{L^2(0,L)}^{2} .
\end{align*}
The sought estimate~\eqref{e:estim-bebete-trace} then follow from~\eqref{e:estim-bebete} and $h\leq 1$.
\enp

We now extend the potentials $\V_{n} ,\V$ as  $\V_{n} ,\V \in C^1_c(( -1, L+1) ; \R )$ (abusing notation slightly) such that $\|\V_{n} - \V\|_{C^1( -1, L+1 )} \to 0$. We define the operator
$$
P_{n}=-h_{n}^{2}\frac{d^{2}}{dx^{2}}+\V_{n} , \quad \text{ acting on } L^2(\R). 
$$
Note that $P_n$ is symmetric on $C^\infty_c(\R)$ since $\V_{n}$ are real-valued.
The equation in~\eqref{e:asspt-psi} together with the jump formula imply that
\bnan
\label{e:eigenP}
P_{n}\udl{\psi_{n}}=-h_{n}^{2}\left(\psi_{n}'(0^{+})\delta_{0}-\psi_{n}'(L^{-})\delta_{L}\right)+ \udl{r_n} , \quad \text{ in } \mathcal{D}'(\R) .
\enan

\begin{corollary}
\label{lmH1trace}
Assume~\eqref{e:asspt-psi}. Then,
\begin{enumerate}
\item \label{i:hn-oscillant} if $r_n= \O_{L^2(0,L)}(1)$ and $\V_n = \O_{L^\infty([0,L])}(1)$, then $h_{n}\big(\udl{\psi_{n}}\big)' = h_{n}\udl{(\psi_{n}')}$ is a bounded sequence in $L^{2}(\R)$ and in particular, 
\begin{align}
\label{e:hosc}
\limsup_{n\to+\infty}\nor{\widehat{\udl{\psi_{n}}}}{L^{2}(|h_{n}\xi|\geq R)}\underset{R\to +\infty}{\longrightarrow}0
\end{align} (where $\hat{u}$ denotes the classical Fourier transform of $u$).
\item \label{i:traces-bounded} if $r_n= \O_{L^2(0,L)}(h_n)$ and $\V_n = \O_{C^1([0,L])}(1)$, then $h_n\psi_n'(0^{+})$ and $h_n\psi_n'(L^{-})$ are bounded sequences in $\R$, and up to a subsequence, there are $\ell_{0}\geq0$ and $\ell_{L}\geq 0$ so that 
\bnan
\label{limtra}
|h_n\psi_n'(0^{+})|^{2}\to \ell_{0}, \quad |h_n\psi_n'(L^{-})|^{2}\to \ell_{L}.
\enan
 Moreover, we have
\bnan
\label{H1oscil}\limsup_{n\to+\infty}\nor{h_n\widehat{ \udl{\psi_{n}'}}}{L^{2}(|h_{n}\xi|\geq R)}\underset{R\to +\infty}{\longrightarrow}0.
\enan
\end{enumerate}
\end{corollary}
Property~\eqref{e:hosc} (resp.~\eqref{H1oscil}) says that the sequence $\udl{\psi_{n}}$ (resp. $h_n\udl{\psi_{n}'}$) is $h_n-$oscillating. This means that the scale $h_n$ ``captures the maximal oscillation rate of the sequence''.
\bnp
Using~\eqref{e:estim-bebete} (applied to $\psi_n$) together with the fact in~\eqref{e:asspt-psi} that $\psi_{n}$ is normalized in $L^{2}$, and the assumption  $r_n= \O_{L^2(0,L)}(1)$, we obtain that $h_{n}\psi_{n}'$ is bounded in $L^2(0,L)$, whence the first statement since $\big(\udl{\psi_{n}}\big)' = \udl{(\psi_{n}')}$ thanks to the Dirichlet boundary condition.
The Plancherel formula then implies that $$\nor{\widehat{\udl{\psi_{n}}}}{L^{2}(|h_{n}\xi|\geq R)}\leq (2\pi)^{-1} R^{-1}\nor{h_{n}\udl{\psi_n}'}{L^{2}}\leq CR^{-1} \to_{R \to + \infty} 0. $$

The fact that $h_n\psi_n'(0^{+})$ and $h_n\psi_n'(L^{-})$ are bounded directly follows from~\eqref{e:estim-bebete-trace} together with the fact that $h_n^{-1}\nor{r_{n}}{L^{2}(0,L)}$ and $\nor{\V_n}{C^1([0,L])}$ are bounded and $\psi_n$ is normalized.

We finally consider the oscillation property for the sequence $h_n \udl{\psi_{n}'}$. Taking $\alpha \in (1/2,1)$, using Equation~\eqref{e:eigenP} and the Plancherel formula, we obtain
\begin{align*}
\nor{\widehat{h_n \udl{\psi_{n}'}}}{L^{2}(|h_{n}\xi|\geq R)}& \leq R^{-1+\alpha}\nor{|h_n\xi|^{-\alpha}h_n^2 \widehat{(\udl{\psi_{n}})''}}{L^{2}(|h_{n}\xi|\geq R)}\\
&\leq 2R^{-1+\alpha}h_{n}^{2}\left(|\psi_{n}'(0^{+})|+|\psi_{n}'(L^{-})|\right)\nor{|h_n\xi|^{-\alpha}\widehat{\delta_0}}{L^{2}(|h_{n}\xi|\geq R)}\\
& \quad + (2\pi)^{-1} R^{-1}\nor{\V_{n}\udl{\psi_{n}}-r_n }{L^{2}}. 
\end{align*}
Since 
$$\nor{|h_n\xi|^{-\alpha}\widehat{\delta_0}}{L^{2}(|h_{n}\xi|\geq R)}= \left(\int_{|h_{n}\xi|\geq R)} |h_n\xi|^{-2\alpha} d\xi\right)^{1/2}  = h_n^{-1/2}\sqrt{2}\left(\int_{R}^\infty |\eta|^{-2\alpha} d\eta\right)^{1/2}  = C_\alpha h_n^{-1/2}R^{-\alpha+1/2} ,$$ 
we then deduce~\eqref{H1oscil} from the facts that $ h_{n}\left(|\psi_{n}'(0^{+})|+|\psi_{n}'(L^{-})|\right)$ and $\nor{\V_{n}\udl{\psi_{n}} -r_n }{L^{2}}$ are bounded.
\enp

\subsection{Localization in the characteristic set}
\label{s:loc-char-set}
The existence of semiclassical measures $\mu$ associated to $(\udl{\psi_n},h_n)_{n\in \N}$ as in~\eqref{e:mu-exists} is classical, see e.g.~\cite[Theorem 5.2]{Zworski:book}.
In this section, we explain how the fact that $\psi_n$ solves Equation~\eqref{e:asspt-psi} (or rather $\udl{\psi_n}$ solves~\eqref{e:eigenP}) relates the associated limit measures $\mu$  to the classical hamiltonian 
$$
p(x,\xi) = \xi^2 + \V(x) .
$$
Remark that in case $\V \in C^\infty(\R)$ and $\V_n=\V$, the function $p(x,\xi)$ is the semiclassical principal symbol of the operator $P_n$. We also denote by $H_p(x,\xi): = 2\xi \d_x - \V'(x)\d_\xi$ the Hamiltonian flow of $p$.
Localization and flow invariance properties for the measures $\mu$ {\em away from the boundary} are proved e.g. in \cite[Theorem 5.5]{Zworski:book} assuming $C^\infty$ regularity. Limited regularity is considered in~\cite{Burq:97}.
Here, we precise these proofs in the case of Dirichlet boundary condition and of family of potentials converging in $C^1$ regularity.

\begin{lemma}
\label{lmsuppmu0}
Assume~\eqref{e:asspt-psi} with $r_n= \O_{L^2(0,L)}(1)$ and $\V_n = \O_{L^\infty([0,L])}(1)$.
Then, the measure $\mu$ in~\eqref{e:mu-exists} is a probability measure supported in the set $[0,L]\times \R_{\xi}$.
\end{lemma}
\bnp
To prove that $\mu$ is a probability measure, we take $\chi ,\chi_L \in C^\infty_c(\R; [0,1])$ such that $\chi =1$ in a neighborhood of $0$, and $\chi_L=1$ in a neighborhood of $[0,L]$, and write (using $\supp(\udl{\psi_n}) \subset [0,L]$)
$$
1 = \nor{\udl{\psi_n}}{L^2(\R)} =  \nor{\chi_L \udl{\psi_n}}{L^2(\R)} \leq \nor{\chi(h_nD/R)\chi_L \udl{\psi_n}}{L^2(\R)} + \nor{\left(1- \chi(h_nD/R)\right)\udl{\psi_n}}{L^2(\R)} .
$$ 
Using Item~\ref{i:hn-oscillant} Corollary~\ref{lmH1trace}, we have $\limsup_{n\to+\infty}\nor{\left(1- \chi(h_nD/R)\right)\udl{\psi_n}}{L^2(\R)}\underset{R\to +\infty}{\longrightarrow}0$, and pseudodifferential calculus yields
$$
\nor{\chi(h_nD/R)\chi_L \udl{\psi_n}}{L^2(\R)}^2 \to_{n\to+\infty} \langle \mu , \chi_L^2 \otimes \chi^2(\cdot /R)\rangle .
$$
We deduce from the above two lines that
$$
1  \leq \langle \mu , \chi_L^2 \otimes \chi^2(\cdot /R)\rangle + o_{R\to \infty}(1) ,
$$
and hence $1 \leq  \langle \mu , \chi_L^2 \otimes 1\rangle \leq 1$ by dominated convergence. This proves both that $\mu$ is a probability measure, and that $\supp (\mu) \subset [0,L]\times \R_{\xi}$.
\enp

\begin{lemma}
\label{lmsuppmu}
Assume~\eqref{e:asspt-psi} with $r_n= \O_{L^2(0,L)}(h_n)$, $\V_n = \O_{C^1([0,L])}(1)$ and $\|\V_{n} - \V\|_{C^0( -1, L+1 )} \to 0$. Then, the measure $\mu$ in~\eqref{e:mu-exists} is a probability measure supported in the set $\{p(x,\xi)=0\}\cap [0,L]\times \R_{\xi}$.
Moreover, for all $a \in C^\infty_c(\R^2;\R)$ such that $a =1$ in neighborhood of $\{p(x,\xi)=0\}\cap [0,L]\times \R_{\xi}$, we have $\nor{\Op_{h_n}(1-a)\udl{\psi_n}}{L^2(\R)}\to 0$ as $n\to +\infty$.
\end{lemma}

Note that the compactness of the set $\{p(x,\xi)=0\}\cap [0,L]\times \R_{\xi} \subset [0,L]\times [-A,A]$, with $A = \sqrt{-\min_{[0,L]}\V}$ thus implies that $\mu \in \E'(\R^2)$, i.e., has compact support.

Note also that the assumption that $r_n= \O_{L^2(0,L)}(h_n)$ can be weakened to $r_n= \O_{L^2(0,L)}(h_n^{1/2+\eps})$ for any $\eps>0$ for the same proof to work (using directly \eqref{e:estim-bebete-trace} instead of Corollary~\ref{lmH1trace} Item~\ref{i:traces-bounded}). We did not try to optimize the proofs in this respect.

\bnp
Let $a\in C^{\infty}_{c}(\R_{x}\times \R_{\xi})$. 
Applying $A=\Op_{h_n}(a)$ to equation \eqref{e:eigenP} and taking the inner product with $\udl{\psi_n} $, we obtain (after having noticed that $\Op_{h_n}(a)$ is a smoothing operator),
\bna
\left(  A P_{n} \udl{\psi_n} ,  \udl{\psi_n} \right)_{L^2(\R)}=-h_{n}^{2}\left( A \left(\psi_{n}'(0^{+})\delta_{0}-\psi_{n}'(L^{-})\delta_{L}\right) ,\udl{\psi_n} \right)_{L^2(\R)}+\petito{1}.
\ena
Corollary~\ref{lmH1trace} Item~\ref{i:traces-bounded} (i.e. boundedness of $h_{n}|\psi_{n}'(0)|$) and continuity of the trace $H^{1/2+\e}(\R) \to \C, u\mapsto u(0)$, gives 
 \begin{align*}
 h_{n}^{2}\left|\left( A(\psi_{n}'(0)\delta_{0}) ,\udl{\psi_n} \right)_{L^2}\right| & =  h_{n}^{2}|\psi_{n}'(0)| \left|\left( A(\delta_{0}) ,\udl{\psi_n} \right)_{L^2}\right| 
 =  h_{n}^{2}|\psi_{n}'(0)| \left|\left< \delta_{0} ,\trans{A} \ovl{\udl{\psi_n}}\right>_{\calS'(\R), \calS(\R)}\right| 
 \\
& \leq C h_{n}\left| ( A^*\udl{\psi_n} )(0)\right|\leq C_\eps h_{n}\nor{A^*\udl{\psi_n}}{H^{1/2+\e}(\R)}\leq C_\eps h_{n}\nor{\udl{\psi_n}}{H^{1/2+\e}(\R)},
 \end{align*} 
 after having used uniform boundedness of $A^*$ on $H^s(\R)$ (classical Sobolev spaces).
The last term is of order $\O_\eps(h_n^{1/2-\eps})$ by interpolation in Corollary~\ref{lmH1trace} between $L^{2}$ and $H^{1}$, and hence converges to zero for $\eps<1/2$.
The same convergence to zero holds for $h_{n}^{2}\left|\left( A(\psi_{n}'(L)\delta_{L}) ,\udl{\psi_n} \right)_{L^2}\right|$, and we have thus proved that for all $a\in C^\infty_c(\R^2)$,  
$\left(  A P_{n}\udl{\psi_n} , \udl{\psi_n} \right)_{L^2(\R)}\to 0$.

For $\epsilon>0$, let $\rho_\epsilon(x) = \frac{1}{\epsilon}\rho(x/\epsilon)$ be an approximation of identity ($\rho \in C^\infty_c(\R)$, $\rho\geq 0$, $\int_\R \rho = 1$). We define $\V^\epsilon := \rho_\epsilon * \V$ and $\V_n^\epsilon := \rho_\epsilon * \V_n$. We notice that for any $\epsilon>0$, we have (under the assumptions of the lemma) that $\V_n^\epsilon = \O_{\epsilon, C^1([0,L])}(1)$ and $\|\V_{n}^\epsilon - \V^\epsilon \|_{C^0( -1, L+1 )} \to 0$ as $n\to +\infty$. Moreover,  $\|\V^\epsilon - \V \|_{C^0( -1, L+1 )} \to 0$ as $\epsilon \to 0$.
We now write 
\begin{align}
\label{e:esti-eps}
\left(  A (h_n^2D_x^2 +\V_n^\epsilon )\udl{\psi_n} , \udl{\psi_n} \right)_{L^2(\R)} = \left(  A P_{n}\udl{\psi_n} , \udl{\psi_n} \right)_{L^2(\R)} +\left(  A (\V_n^\epsilon -\V_n) \udl{\psi_n} , \udl{\psi_n}\right)_{L^2(\R)} .
\end{align} 
The first term in the right hand-side converges to zero, whereas the second term is bounded by $$\nor{A}{\L(L^2)}\nor{\V_n^\epsilon -\V_n}{C^0}\to \nor{A}{\L(L^2)}\nor{\V^\epsilon -\V}{C^0}, \quad \text{ as } n\to + \infty.
$$ Pseudodifferential calculus (composition rule) in the left hand-side of~\eqref{e:esti-eps}, together with the fact that $\|\V_{n}^\epsilon - \V^\epsilon \|_{C^0( -1, L+1 )} \to 0$ as $n\to +\infty$ and the definition of $\mu$ imply that it converges towards $\langle \mu , (|\xi|^2 + \V^\epsilon)a \rangle$. We have thus obtained that
$$
 \left| \langle \mu , (|\xi|^2 + \V^\epsilon)a \rangle \right| \leq \nor{A}{\L(L^2)}\nor{\V^\epsilon -\V}{C^0} \to 0 , \text{ as }\epsilon \to 0^+.
$$
Since  $\langle \mu , (|\xi|^2 + \V^\epsilon)a \rangle\underset{\epsilon \to 0^+} { \rightarrow}\langle \mu , (|\xi|^2 + \V)a \rangle$, we have obtained $\langle \mu , pa \rangle = 0$ for all $a \in  C^\infty_c(\R^2)$. This implies that $\supp(\mu) \subset p^{-1}(\{0\})$, and concludes the proof of first statement of the lemma. 
 
 Concerning the second statement, using pseudodifferential calculus and the normalization of the $\psi_n$'s, we have 
 \begin{align*}
 \nor{\Op_{h_n}(1-a)\udl{\psi_n}}{L^2(\R)}^2 
& =\left( \Op_{h_n}\left((1-a)^2 \right)\udl{\psi_n},  \udl{\psi_n} \right)_{L^2(\R)}  + \O(h_n) \\
 & =\nor{ \udl{\psi_n}}{L^2(\R)}^2 + \left( \Op_{h_n}(-2a +a^2)\udl{\psi_n},  \udl{\psi_n} \right)_{L^2(\R)}  + \O(h_n) \\
 & \to 1 + \langle \mu , -2a +a^2 \rangle .
 \end{align*}
 Recalling that $a=1$ in a neighborhood of $\supp(\mu)$ and that $\mu$ is a probability measure, we have $ \langle \mu , -2a +a^2 \rangle = \langle \mu , -1\rangle = -1$, whence the sought result.
\enp

\subsection{Propagation of the measure}
\label{s:propagation}
We next want to investigate propagation properties for the measure $\mu$, and start with a propagation statement ``away from the boundary''.
\begin{lemma}
\label{lmsupHpmu}
Under the assumptions of Lemma \ref{lmsuppmu}, with $\V \in C^1_c([-1,L+1])$, the distribution $H_{p}\mu$ defined by
\begin{align}
\label{e:def-Hpmu}
\langle H_{p}\mu, a \rangle_{\mathcal{D}'(\R^2),\mathcal{D}(\R^2)}  : = - \langle \mu , (2\xi \d_x- \V' \d_{\xi}) a \rangle_{\mathcal{M}(\R^2),C^0(\R^2)} , \quad a \in C^\infty_c(\R^2) ,
\end{align}
 is of order at most $1$.
 If moreover $r_n= o_{L^2(0,L)}(h_n)$ and $\|\V_{n} - \V\|_{C^1( -1, L+1 )} \to 0$, then 
\begin{align}
\label{e:supp-hpmu-intern}
\supp( H_{p}\mu ) \subset \{\xi^{2}+\V(x)=0\}\cap \big( \{0,L\}\times \R_{\xi} \big).
\end{align}
\end{lemma}
The support statement~\eqref{e:supp-hpmu-intern} in Lemma~\ref{lmsupHpmu} says that the measure $\mu$ is $H_p$-invariant ``away from the boundary'' of the interval $[0,L]$.
The proof is classical in case $\V_n=\V$ is smooth, but requires some care in the present limited regularity setting. 
\bnp
Since $\V\in C^1(\R)$, we have from the definition~\eqref{e:def-Hpmu} that  $|\langle H_{p}\mu, a \rangle_{\mathcal{D}'(\R^2),\mathcal{D}(\R^2)} | \leq  C_K \| a\|_{C^1(\R^2)}$ for all $a \in C^\infty_c(K)$, $K \subset \R^2$ compact. Hence $H_{p}\mu$ is a distribution of order $1$. 

Let us now turn to the support property~\eqref{e:supp-hpmu-intern}. Lemma \ref{lmsuppmu} first implies that  $H_{p}\mu$ is supported in the set $\{\xi^{2}+\V(x)=0\}\cap \big( [0,L]\times \R_{\xi} \big)$. It is therefore sufficient to prove that 
\begin{align}
\label{e:remains-to-be-proved}
\langle H_{p}\mu, a \rangle_{\mathcal{D}'((0,L)\times \R_{\xi}),\mathcal{D}((0,L)\times \R_{\xi})}=0
\end{align}  for any $a \in C^\infty_c((0,L)\times \R_{\xi})$ (with $H_p\mu$ defined in the sense of~\eqref{e:def-Hpmu}). By density of vector space spanned by tensor products of smooth functions in $C^1_c((0,L)\times \R_{\xi})$, it is enough to prove~\eqref{e:remains-to-be-proved} for test functions $a$ under the form $a(x,\xi)=\chi_{1}(x)\chi_{2}(\xi)$ with $\chi_{1}\in C^\infty_c(0,L)$ and $\chi_{2}\in C^\infty_c( \R_{\xi})$.

As in the proof of Lemma~\ref{lmsuppmu}, for $\epsilon>0$, we let $\rho_\epsilon(x) = \frac{1}{\epsilon}\rho(x/\epsilon)$ be an approximation of identity and define $\V^\epsilon := \rho_\epsilon * \V$ and $\V_n^\epsilon := \rho_\epsilon * \V_n$. We notice that for any $\epsilon>0$, we have $\|\V_{n}^\epsilon - \V^\epsilon \|_{C^1( -1, L+1 )} \to 0$ as $n\to + \infty$, and $\|\V^\epsilon - \V \|_{C^1( -1, L+1 )} \to 0$ as $\epsilon \to 0$.
We also set $A=\chi_{2}(h_{n}D_{x})\chi_{1}(x)$.
 The proof of~\eqref{e:remains-to-be-proved} consists in computing in two different ways the limit of the quantity 
\begin{align}
\label{e:LA-def0interne}
L_{A}(h_n) :=\frac{1}{h_{n}}\left< AP_{n}\udl{\psi_n} , \ovl{\udl{\psi_n}}  \right>_{\mathcal{S}(\R),\mathcal{S'}(\R)}-\frac{1}{h_{n}}\left<  A \udl{\psi_n}  , P_{n}  \ovl{\udl{\psi_n}} \right>_{\mathcal{S}(\R),\mathcal{S'}(\R)} ,
\end{align}
which makes sense since $A \udl{\psi_n}\in \mathcal{S}(\R)$ and $AP_{n}\udl{\psi_n}\in \mathcal{S}(\R)$.
 Using that $P_{n}$ is formally selfadjoint together with pseudodifferential rules, we have on the one hand
\begin{align}
\label{e:LA-definterne}
L_{A}(h_n)  =\frac{1}{h_{n}}\left( [A,P_{n}]\udl{\psi_n} , \udl{\psi_n} \right)_{L^2}
= \frac{1}{h_{n}}\left( [A,P_{n}^\epsilon]\udl{\psi_n} , \udl{\psi_n} \right)_{L^2} + \frac{1}{h_{n}}\left( [A,P_{n} - P_{n}^\epsilon]\udl{\psi_n} , \udl{\psi_n} \right)_{L^2} , 
\end{align}
where $P_{n}^\epsilon=-h_{n}^{2}\frac{d^{2}}{dx^{2}}+\V_{n}^\epsilon$.
We first study the first term in~\eqref{e:LA-definterne}. For fixed $\epsilon>0$, $[A,P_{n}^\epsilon]$ is a semiclassical operator, so we can write
\begin{align}
\label{e:LA-bisinterne}
\frac{1}{h_{n}}\left( [A,P_{n}^\epsilon]\udl{\psi_n} , \udl{\psi_n} \right)_{L^2} \underset{n\to+\infty}{\to} \left<\mu,\frac1{i} \{a ,p^\epsilon\}  \right> 
= -\frac{1}{i}\left<\mu,H_{p^\epsilon} a \right> \to_{\epsilon \to 0^+} -\frac{1}{i}\left<\mu,H_{p} a \right> .
\end{align}
Concerning the second term in~\eqref{e:LA-definterne}, we have
$\frac{1}{h_{n}}\left( [A,P_{n} - P_{n}^\epsilon]\udl{\psi_n} , \udl{\psi_n} \right)_{L^2}
= \frac{1}{h_{n}}\left( [A,\V_{n} - \V_{n}^\epsilon]\udl{\psi_n} , \udl{\psi_n} \right)_{L^2}$
where, using the product form of $A$,   
\begin{align*}
[A,\V_{n} - \V_{n}^\epsilon] & =  \chi_{2}(h_{n}D_{x})[ \chi_{1}(x),\V_{n} - \V_{n}^\epsilon]
+[\chi_{2}(h_{n}D_{x}),\V_{n} - \V_{n}^\epsilon]\chi_{1}(x) =[\chi_{2}(h_{n}D_{x}),\V_{n} - \V_{n}^\epsilon]\chi_{1}(x).
\end{align*}
As a consequence, recalling that $\udl{\psi_n}$ is normalized in $L^2(\R)$, we obtain
\begin{align}
\label{e:LA-estiinterne}
\left| \frac{1}{h_{n}}\left( [A,P_{n} - P_{n}^\epsilon]\udl{\psi_n} , \udl{\psi_n} \right)_{L^2}\right|
& \leq \frac{1}{h_{n}}\nor{[\chi_{2}(h_{n}D_{x}),\V_{n} - \V_{n}^\epsilon]}{\L(L^2)} \nonumber\\
& \leq C \nor{\d_x(\V_{n} - \V_{n}^\epsilon)}{L^\infty} \underset{n\to+\infty}{\to} C \nor{\d_x(\V - \V^\epsilon)}{L^\infty}  \to_{\epsilon \to 0^+} 0,
\end{align}
 where we used Lemma~\ref{l:comm-pas-reg}  below with $\eps = h_n$. Combining~\eqref{e:LA-definterne} with~\eqref{e:LA-bisinterne} and \eqref{e:LA-estiinterne}, and letting $n \to+\infty$ and then $\epsilon \to 0^{+}$ (recall that $L_{A}(h_n)$ is independent of $\epsilon$), we have obtained for $a(x,\xi)=\chi_{1}(x)\chi_{2}(\xi)$ with $\chi_{1}\in C^\infty_c(0,L)$ and $\chi_{2}\in C^\infty_c( \R_{\xi})$,
\begin{align}
\label{e:tototo-tututu}
L_{A}(h_n) \underset{n\to+\infty}{\to} -\frac{1}{i}\left<\mu,H_{p} a \right> .
\end{align}
We now compute $L_{A}(h_n)$ in~\eqref{e:LA-def0interne} using Equation~\eqref{e:eigenP}. Using moreover that $A=\chi_{2}(h_{n}D_{x})\chi_{1}(x)$ with $\supp (\chi_1) \subset (0,L)$, together with the pseudolocality of $A^{*}$, we obtain
\begin{align*}
L_{A}(h_n)=\frac{1}{h_{n}}\left( A \udl{r_n} , \udl{\psi_n} \right)_{L^2(\R)}-\frac{1}{h_{n}}\left(  \udl{\psi_n}  ,A^* \udl{r_n} \right)_{L^2(\R)} + \mathcal{O}(h_n^\infty)\| \udl{\psi_n}\|_{L^2(\R)}(|\psi_n'(0)|  +|\psi_n'(L)| ) .
\end{align*}
Then, $L^2$ normalization of $\udl{\psi_n}$ together with $L^2$ boundedness of $A$ and the assumption $r_n= o_{L^2(0,L)}(h_n)$ imply that the first two terms converge to zero as $n \to +\infty$. 
Item~\ref{i:traces-bounded} of Lemma~\ref{lmH1trace} implies that the last term converges to zero as well.
Combined with~\eqref{e:tototo-tututu}, we have thus obtained~\eqref{e:remains-to-be-proved} for all $a$ in product form, and finally for all  $a \in C^\infty_c((0,L)\times \R_{\xi})$ by density. This concludes the proof of~\eqref{e:supp-hpmu-intern}, and thus of the lemma.
\enp

As a preliminary for propagation properties, we first prove that the convergence in~\eqref{e:mu-exists} holds not only for compactly supported symbols $a$ but also for symbols of order $2$. For $m\in \R$, we shall say that $a \in S^m(\R^2)$ if $|\d_x^\alpha\d_\xi^\beta a(x,\xi, h) | \leq C_{\alpha,\beta}\langle \xi \rangle^{m-\beta}$ for all $\alpha, \beta \in \N, (x,\xi) \in \R^2, h \in (0,1]$ (note that such symbols depend implicitly on $h$, with uniform bounds).

\begin{lemma}
\label{l:conv-S1}
Assume~\eqref{e:asspt-psi} with $r_n= \O_{L^2(0,L)}(h_n)$, $\V_n = \O_{C^1([0,L])}(1)$ and $\|\V_{n} - \V\|_{C^0( -1, L+1 )} \to 0$.
Then, for all $a \in S^2(T^*\R)$ independent of $h$, we have
\begin{align}
\label{e:H-1H1A}
\left< \Op_{h_n}(a) \udl{\psi_n} ,  \udl{\psi_n} \right>_{H^{-1},H^1} \to \langle \mu , a \rangle_{\E'(\R^2),\E(\R^2)}  . 
\end{align}
\end{lemma}
We denote for $s \in \R$  by $$
\|u\|_{H^s_h}^2  = \int_\R (1+h^2|\xi|^2)^{s/2}|\hat{u}(\xi)|^2 d\xi , \quad \text{ where } \hat{u}(\xi) = \int_{\R} e^{-i x \xi} u(x) dx , \quad u \in \mathcal{S}(\R) ,
$$ the usual semiclassical Sobolev norm.
Note that in Expression~\eqref{e:H-1H1A}, $\Op_{h_n}(a) \udl{\psi_n} \in H^{-1}$, and is bounded uniformly in $H^{-1}_{h_n}$ since $\udl{\psi_n}$ and $h_n \udl{\psi_n}'$ are bounded in $L^2(\R)$, see Corollary~\ref{lmH1trace}.
In particular, one can replace in~\eqref{e:H-1H1A} $a \in S^2(T^*\R)$ independent of $h$ by $a + \eps(h) b$ with $b \in S^2(T^*\R)$ possibly depending on $h$ (with uniformly bounded seminorms in this class) and $\eps(h) \to 0$.

Here and below, we take the convention that duality brackets between $H^{-1}$ and $H^1$ or between $H^{-1}_h$ and $H^1_h$ in~\eqref{e:H-1H1A} are $\C-$linear in the first variable and $\C-$antilinear in the second variable.

\medskip
Before giving the proof of Lemma~\ref{l:conv-S1}, we give the following corollary which is actually the last item of Proposition~\ref{propmeasureapp} (and is valid under less restrictive assumptions).
\begin{corollary}
\label{c:proj-mes}
We have 
$|\udl{\psi_n}(x)|^2 dx \rightharpoonup \mathfrak{m}$ where the nonnegative Radon measure $\mathfrak{m}$ on $\R$ is given by $\mathfrak{m}=\pi^{*}\mu$.
\end{corollary}
\bnp[Proof of Corollary~\ref{c:proj-mes} from Lemma~\ref{l:conv-S1}]
For any $\varphi\in C^{\infty}_{c}(\R)$, we can apply Lemma \ref{l:conv-S1} to $a=\varphi\circ \pi \in S^0(T^*\R)$, to obtain $\int_{\R}\varphi(x)|\udl{\psi_n}(x) |^{2} ~dx \to \langle \mu , \varphi \circ \pi \rangle_{\E'(\R^2),\E(\R^2)}=\langle \mathfrak{m} , \varphi  \rangle_{\E'(\R),\E(\R)}$ by definition of $\mathfrak{m}$. By density of $C^{\infty}_{c}(\R)$ in $C^{0}_{c}(\R)$, this implies the result.
\enp

\bnp[Proof of Lemma \ref{l:conv-S1}]
We choose $\phi \in C^\infty_c(\R^2)$ real-valued such that $\phi=1$ in a neighborhood of $\supp(\mu)$.
We decompose 
$$
\left< \Op_{h_n}(a) \udl{\psi_n} ,  \udl{\psi_n} \right>_{H^{-1},H^1} = 
\left<\Op_{h_n}(\phi)\Op_{h_n}(a)  \udl{\psi_n} ,  \udl{\psi_n} \right>_{H^{-1},H^1} +
\left<(1-\Op_{h_n}(\phi))  \Op_{h_n}(a) \udl{\psi_n} , \udl{\psi_n} \right>_{H^{-1},H^1} .
$$
We first notice that we have on the one hand (for any $a \in S^m(\R^2)$ with principal part independent of $h$)
\begin{align*}
\left(  \Op_{h_n}(\phi) \Op_{h_n}(a)  \udl{\psi_n} ,  \udl{\psi_n} \right)_{L^2}&  = \left( \Op_{h_n}(a\phi)  \udl{\psi_n} ,  \udl{\psi_n} \right)_{L^2} + \O(h_n)  \underset{n\to+\infty}{\to} \left<\mu,\phi a  \right>  = \left<\mu, a  \right> , 
\end{align*}
using pseudodifferential calculus and the support properties of $\phi$. 

To conclude the proof, it suffices to prove
\begin{align}
\label{e:target}
\left<(1-\Op_{h_n}(\phi))  \Op_{h_n}(a) \udl{\psi_n} , \udl{\psi_n} \right>_{H^{-1},H^1} \to 0 .
\end{align}
We first prove the intermediate statement that 
\begin{align}
\label{e:1-chi-0}
\nor{(1-\Op_{h_n}(\phi)) h_{n}\udl{\psi_n'} }{L^2} \underset{n\to+\infty}{\to}0 .
\end{align} 
To prove~\eqref{e:1-chi-0}, we decompose for $\chi \in C^\infty_c(\R)$ equal to one near zero,
$$
(1-\Op_{h_n}(\phi)) h_{n}\udl{\psi_n'}=(1-\Op_{h_n}(\phi)) \chi(R^{-1}h_{n}D_{x})h_{n}\udl{\psi_n'}+(1-\Op_{h_n}(\phi))(1- \chi(R^{-1}h_{n}D_{x}))h_{n}\udl{\psi_n'} ,
$$
 for $R$ large. For the second term, we have
\bna
\limsup_{n\to+\infty}\nor{(1-\Op_{h_n}(\phi))(1- \chi(R^{-1}h_{n}D_{x}))h_{n}\udl{\psi_n'}}{L^{2}}
\leq C \limsup_{n\to+\infty}\nor{(1- \chi(R^{-1}h_{n}D_{x}))h_{n}\udl{\psi_n'}}{L^{2}}  \underset{R\to+\infty}{\to}0,
\ena
using \eqref{H1oscil}. As for the first term, using that $\udl{\psi_n} $ is supported in $[0,L]$, we write for any $R>0$
\bna
\nor{(1-\Op_{h_n}(\phi)) \chi(R^{-1}h_{n}D_{x})h_{n}\udl{\psi_n'}}{L^{2}}= \left(B_R \udl{\psi_n} ,  \udl{\psi_n} \right)_{L^2} 
\ena
for $B_R=\chi(x/L)\chi(R^{-1}h_{n}D_{x})(1-\Op_{h_n}(\phi)^{*}) (1-\Op_{h_n}(\phi)) \chi(R^{-1}h_{n}D_{x})$, that is a semiclassical pseudodifferential operator in $\Op_{h_n}(S^0)$. 
Writing $B_R = \Op_{h}(b_R) + h \Op_{h}(S^{-1})$ with $b_R (x,\xi)= \chi(x/L)\chi(R^{-1}\xi)^2 (1-\phi(x,\xi))^2$, we have obtained
$$
 \left(B_R \udl{\psi_n} ,  \udl{\psi_n} \right)_{L^2}  \to \langle \mu , b_R \rangle = 0 , 
$$
since $\phi=1$ in a neighborhood of $\supp(\mu)$, for any $R>0$. Combining the above three lines, we have proved~\eqref{e:1-chi-0}.

We finally prove~\eqref{e:target} for $a \in S^2$. Denoting by $\transp\Op_{h_n}(\phi)$ the transpose of $\Op_{h_n}(\phi)$ for the duality bracket between $H^{-1}_h$ and $H^1_h$, we have
\begin{align*}
\left|\left<(1-\Op_{h_n}(\phi))  \Op_{h_n}(a) \udl{\psi_n} , \udl{\psi_n} \right>_{H^{-1},H^1} \right|
& = \left| \left< \Op_{h_n}(a) \udl{\psi_n} , (1- \transp\Op_{h_n}(\phi)) \udl{\psi_n} \right>_{H^{-1}_h, H^{1}_h} \right| \\
&\leq \left\|  \Op_{h_n}(a) \udl{\psi_n} \right\|_{H^{-1}_h} 
 \left\| (1-\transp\Op_{h_n}(\phi)) \udl{\psi_n}\right\|_{H^1_h}   \\
&\leq C \left\| \udl{\psi_n} \right\|_{H^{1}_h} \left(  \nor{(1-\transp\Op_{h_n}(\phi)) h_{n}\udl{\psi_n}' }{L^2} + \nor{(1-\transp\Op_{h_n}(\phi)) \udl{\psi_n} }{L^2} \right) \\
&  \underset{n\to+\infty}{\to}0 ,
\end{align*}
where we have used~\eqref{e:1-chi-0}, the fact that $\Op_{h}(a) : H^1_{h} \to H^{-1}_{h}$ uniformly in $h$, and Lemma~\ref{lmsuppmu} for the last convergence. This concludes the proof of the lemma. 
\enp

Note that the right hand-side of~\eqref{e:H-1H1A} makes sense for any $a \in C^\infty(\R^2)$, using that $\mu$ is compactly supported. Convergence in~\eqref{e:H-1H1A} however uses $a \in S^2(\R^2)$.
\begin{lemma}
\label{lmderivbord}
Assume~\eqref{e:asspt-psi} with $r_n= o_{L^2(0,L)}(h_n)$ and $\|\V_{n} - \V\|_{C^1( -1, L+1 )} \to 0$. Then, for all $a_{0}$, $a_{1} \in C^{\infty}_{0}(\R_{x})$ real valued and $a(x,\xi)=a_{0}(x)+ a_{1}(x)\xi$, the measure $\mu$ in~\eqref{e:mu-exists} satisfies
\bna
\left<\mu,H_{p} a\right>= -\ell_{0}a_{1}(0)+\ell_{L}a_{1}(L).
\ena
\end{lemma}
In the next proofs we need an  function $\chi$ such that 
\begin{equation}
\label{e:def-chi}\chi\in C^{\infty}_{c}((-2,2) ;[0,1]), \quad \chi \text{ even}, \quad  \chi(s)=1 \text{ for }|s|\leq 1 ,
\end{equation} 
When $\e>0$, is given, we will denote the function $\chi_{\e}(s)=\chi(\e s)$.

The proof of Lemma~\ref{lmderivbord} follows the general scheme of that of Lemma~\ref{lmsupHpmu}, but we now need an extra care to handle the boundary terms.
\bnp
We set $A=\chi(h_{n}^{3}D_{x})A_0$ with $A_0 :=  a_{0}(x)+ a_{1}(x)h_{n}D_{x}$.
As in the proofs of Lemmata~\ref{lmsuppmu} and~\ref{lmsupHpmu}, for $\epsilon>0$, we let $\rho_\epsilon(x) = \frac{1}{\epsilon}\rho(x/\epsilon)$ be an approximation of identity and define $\V^\epsilon := \rho_\epsilon * \V$ and $\V_n^\epsilon := \rho_\epsilon * \V_n$. We notice that for any $\epsilon>0$, we have $\|\V_{n}^\epsilon - \V^\epsilon \|_{C^1( -1, L+1 )} \to 0$ as $n\to + \infty$, and $\|\V^\epsilon - \V \|_{C^1( -1, L+1 )} \to 0$ as $\epsilon \to 0$.

 The proof consists in computing in two different ways the limit of the quantity 
\begin{align}
\label{e:LA-def0}
L_{A}(h_n) :=\frac{1}{h_{n}}\left< AP_{n}\udl{\psi_n} , \ovl{\udl{\psi_n}}  \right>_{\mathcal{S}(\R),\mathcal{S'}(\R)}-\frac{1}{h_{n}}\left<  A \udl{\psi_n}  , P_{n}  \ovl{\udl{\psi_n}} \right>_{\mathcal{S}(\R),\mathcal{S'}(\R)} ,
\end{align}
which makes sense since $A \udl{\psi_n}\in \mathcal{S}(\R)$ and $AP_{n}\udl{\psi_n}\in \mathcal{S}(\R)$.
 Using that $P_{n}$ is formally selfadjoint together with pseudodifferential rules, we have on the one hand
\begin{align}
\label{e:LA-def}
L_{A}(h_n)  =\frac{1}{h_{n}}\left( [A,P_{n}]\udl{\psi_n} , \udl{\psi_n} \right)_{L^2}
= \frac{1}{h_{n}}\left( [A,P_{n}^\epsilon]\udl{\psi_n} , \udl{\psi_n} \right)_{L^2} + \frac{1}{h_{n}}\left( [A,P_{n} - P_{n}^\epsilon]\udl{\psi_n} , \udl{\psi_n} \right)_{L^2} , 
\end{align}
where $P_{n}^\epsilon=-h_{n}^{2}\frac{d^{2}}{dx^{2}}+\V_{n}^\epsilon$.

We first study the first term in~\eqref{e:LA-def}. 
Recalling $A=\chi(h_{n}^{3}D_{x})A_0$, we decompose  
$$[A,P_{n}^\epsilon] =  \chi(h_{n}^{3}D_{x}) [A_0,P_{n}^\epsilon] + [\chi(h_{n}^{3}D_{x}),P_{n}^\epsilon]  A_0 .$$ On the one hand, we have $[\chi(h_{n}^{3}D_{x}),P_{n}^\epsilon] = [\chi(h_{n}^{3}D_{x}),\V_{n}^\epsilon] =\O_{\L(L^2)}(h_n^3)$ according to pseudodifferential calculus (or Lemma~\ref{l:comm-pas-reg} below), and thus 
\begin{align*}
\left| \frac{1}{h_{n}}\left( [\chi(h_{n}^{3}D_{x}),P_{n}^\epsilon]  A_0 \udl{\psi_n} , \udl{\psi_n} \right)_{L^2} \right|
\leq C_\epsilon h_n^2 \|A_0 \udl{\psi_n}\|_{L^2}  \|\udl{\psi_n}\|_{L^2} = \O_\epsilon(h_n^2) ,
\end{align*}
according to Corollary~\ref{lmH1trace} and the definition of $A_0$. On the other hand, we have 
\begin{align*}
\frac{1}{h_{n}}\left<\chi(h_{n}^{3}D_{x}) [A_0,P_{n}^\epsilon]\udl{\psi_n} , \udl{\psi_n} \right>_{H^{-1},H^1} 
= \left< \frac{1}{h_{n}}[A_0,P_{n}^\epsilon]\udl{\psi_n} ,  \udl{\psi_n} \right>_{H^{-1},H^1} + R_n^{\epsilon} ,
\end{align*}
with 
\begin{align*}
|R_n^{\epsilon} | = \left|\left< \frac{1}{h_{n}}[A_0,P_{n}^\epsilon]\udl{\psi_n} ,(1-\chi(h_{n}^{3}D_{x}) ) \udl{\psi_n} \right>_{H^{-1},H^1} \right| \leq C_\eps 
\nor{\udl{\psi_n}}{H^1_{h_n}} \nor{(1-\chi(h_{n}^{3}D_{x}) ) \udl{\psi_n} }{H^1_{h_n}},
\end{align*}
using that $\frac{1}{h_{n}}[A_0,P_{n}^\epsilon] \in \Op_{h_n}(S^2)$ (actually, it is a semiclassical differential operator of order $2$, that is finite sum of terms of the form $c_{jk}(x)h_n^k D_x^j$, $k\geq j$, $0\leq j\leq 2$). We conclude that $R_n^{\epsilon}$ converges to zero  as $n\to \infty$  thanks to Corollary~\ref{lmH1trace}.

Combining the above lines and using Lemma~\ref{l:conv-S1} (and the remark thereafter), we have obtained that the first term in~\eqref{e:LA-def} satisfies
\begin{align}
\label{e:LA-bis}
\frac{1}{h_{n}}\left( [A,P_{n}^\epsilon]\udl{\psi_n} , \udl{\psi_n} \right)_{L^2} & = \left< \frac{1}{h_{n}}[A_0,P_{n}^\epsilon]\udl{\psi_n} ,  \udl{\psi_n} \right>_{H^{-1},H^1} +o_\epsilon(1) \nonumber \\
&\underset{n\to+\infty}{\to} \left<\mu,\frac1{i} \{a ,p^\epsilon\}  \right> 
= -\frac{1}{i}\left<\mu,H_{p^\epsilon} a \right> \to_{\epsilon \to 0^+} -\frac{1}{i}\left<\mu,H_{p} a \right> .
\end{align}
Concerning the second term in~\eqref{e:LA-def}, we have
$\frac{1}{h_{n}}\left( [A,P_{n} - P_{n}^\epsilon]\udl{\psi_n} , \udl{\psi_n} \right)_{L^2}
= \frac{1}{h_{n}}\left( [A,\V_{n} - \V_{n}^\epsilon]\udl{\psi_n} , \udl{\psi_n} \right)_{L^2}$
where 
\begin{align*}
[A,\V_{n} - \V_{n}^\epsilon] & =  \chi(h_{n}^{3}D_{x})[\left( a_{0}(x)+ a_{1}(x)h_{n}D_{x}\right),\V_{n} - \V_{n}^\epsilon]
+[\chi(h_{n}^{3}D_{x}),\V_{n} - \V_{n}^\epsilon]\left( a_{0}(x)+ a_{1}(x)h_{n}D_{x}\right) \\
& =  \chi(h_{n}^{3}D_{x})a_{1}(x) \frac{h_{n}}{i}\d_x(\V_{n} - \V_{n}^\epsilon)
+[\chi(h_{n}^{3}D_{x}),\V_{n} - \V_{n}^\epsilon]\left( a_{0}(x)+ a_{1}(x)h_{n}D_{x}\right).
\end{align*}
As a consequence, recalling that $ \udl{\psi_n}$ and $h_n \udl{\psi_n}'$ are bounded in $L^2$, we obtain
\begin{align}
\label{e:LA-esti}
\left| \frac{1}{h_{n}}\left( [A,P_{n} - P_{n}^\epsilon]\udl{\psi_n} , \udl{\psi_n} \right)_{L^2}\right|
& =\left|  \frac{1}{i}\left(a_{1}(x) \d_x(\V_{n} - \V_{n}^\epsilon)
 \udl{\psi_n} ,\chi(h_{n}^{3}D_{x}) \udl{\psi_n} \right)_{L^2} \right. \nonumber\\
 & \quad \left. 
-\frac{1}{h_{n}} \left(\left( a_{0}(x)+ a_{1}(x)h_{n}D_{x}\right) \udl{\psi_n} ,[\chi(h_{n}^{3}D_{x}),\V_{n} - \V_{n}^\epsilon] \udl{\psi_n} \right)_{L^2} \right| \nonumber\\
& \leq C \nor{\d_x(\V_{n} - \V_{n}^\epsilon)}{L^\infty} + \frac{C}{h_{n}}\nor{[\chi(h_{n}^{3}D_{x}),\V_{n} - \V_{n}^\epsilon]}{\L(L^2)}\nonumber \\
& \leq C \nor{\d_x(\V_{n} - \V_{n}^\epsilon)}{L^\infty} +Ch_n^2 \nor{\d_x(\V_{n} - \V_{n}^\epsilon)}{L^\infty}  \underset{n\to+\infty}{\to} C \nor{\d_x(\V - \V^\epsilon)}{L^\infty}  ,
\end{align}
 where we used Lemma~\ref{l:comm-pas-reg}  below with $\eps = h_n^3$ in the last line. Combining now~\eqref{e:LA-def} with~\eqref{e:LA-bis}, \eqref{e:LA-esti}, and the fact that $\nor{\d_x(\V - \V^\epsilon)}{L^\infty} \to 0$ as $\epsilon \to 0$, we have obtained 
\begin{align} 
 \label{e:LA-on-one-side}
 L_{A}(h_n) \underset{n\to+\infty}{\to} -\frac{1}{i}\left<\mu,H_{p} a \right> .
 \end{align}
 
 We now compute $ L_{A}(h_n)$ defined in~\eqref{e:LA-def0} in a different way using the equation~\eqref{e:eigenP}. We obtain
 \begin{align}
 \label{e:bdry-LA}
 L_{A}(h_n)&=-h_{n} \left< A\left(\psi_{n}'(0^{+})\delta_{0}-\psi_{n}'(L^{-})\delta_{L}\right) , \ovl{\udl{\psi_n}} \right>_{\mathcal{S}(\R),\mathcal{S'}(\R)} \nonumber\\
  & \quad +h_{n}\left< A\udl{\psi_n} , \left(\ovl{\psi_{n}}'(0^{+})\delta_{0}-\ovl{\psi_{n}}'(L^{-})\delta_{L}\right) \right>_{\mathcal{S}(\R),\mathcal{S'}(\R)}+\petito{1} \nonumber\\
 &=h_{n} \left[-\psi_{n}'(0^{+}) (\ovl{A^{*}\udl{\psi_{n}}})(0)+(A\udl{\psi_{n}})(0)\overline{\psi_{n}}'(0^+)\right] \nonumber \\
  & \quad +h_{n}\left[\psi_{n}'(L^{-})( \ovl{A^{*}\udl{\psi_n}})(L)-(A \udl{\psi_n})(L)\overline{\psi_{n}}'(L^{-})\right]+\petito{1}.
 \end{align}
 We now only treat the boundary terms at $0$; the boundary terms at $L$ being handled similarly. Recalling the definition of $A$ at the beginning of the proof, we have
 $$
 A^* =\left( a_{0}(x) + \frac{h_n}{i}a_1'(x)+ a_{1}(x)h_{n}D_{x}\right) \chi(h_{n}^{3}D_{x}) .
 $$
 As a consequence, we have 
 \begin{align}
 \label{e:bdry-LA-bis}
 (\ovl{A^{*}\udl{\psi_{n}}})(0)&=\left( a_{0}(0) + \frac{h}{i}a_1'(0) \right)\left[\chi(h_{n}^{3}D_{x})\ovl{\udl{\psi_{n}}}\right](0)+a_{1}(0)\left[\chi(h_{n}^{3}D_{x})h_{n}D_{x}\ovl{\udl{\psi_{n}}}\right](0) , \\
  \label{e:bdry-LA-ter}
  (A\udl{\psi_{n}})(0)&=\left[\chi(h_{n}^{3}D_{x})\left( a_{0}\udl{\psi_n} \right) \right](0)+ \left[\chi(h_{n}^{3}D_{x}) \left(a_{1}h_{n}D_{x}\udl{\psi_{n}}\right)\right](0) . 
   \end{align}
 It is now possible to apply Lemma \ref{lmtracejump} below with $\e=h_{n}^{3}$ with $f=\udl{\psi_{n}}$ or $h_{n} D_{x}\udl{\psi_{n}}$ or $a_{0}(x)\udl{\psi_{n}}$ or $a_{1}(x)h_{n}D_{x}\udl{\psi_{n}}$. For instance, using that $\udl{\psi_n}(0)=0$, we have 
 $$
 \left| \left[\chi(h_{n}^{3}D_{x})\left( a_{0}\udl{\psi_n} \right) \right](0) \right| \leq C h_n^\frac32 \left( \nor{ (a_{0}\udl{\psi_n})' }{L^2}+  \nor{ a_{0}\udl{\psi_n} }{L^2}\right)  \leq C h_n^\frac12 , 
 $$
 on account to Corollary~\ref{lmH1trace}. Similarly, according to  Lemma \ref{lmtracejump}, we have 
 \begin{align*}
 &\left[\chi(h_{n}^{3}D_{x}) \left(a_{1}h_{n}D_{x}\udl{\psi_{n}}\right)\right](0) = \frac12 a_{1}(0)h_{n}D_{x}\psi_{n}(0^+) +s_n ,  \\
  & \text{with } \quad |s_n| \leq C h_n^\frac52 \left( \nor{ (a_{1}\psi_n')' }{L^2(0,L)}+  \nor{a_{1}\psi_n' }{L^2(0,L)}\right)  \leq C h_n^\frac12 , 
 \end{align*}
 where we used the equation~\eqref{e:asspt-psi}.
 Note that the power $3$ in $\chi(h_{n}^{3}D_{x})$ has been chosen so that to handle the remainder terms. 
Collecting all terms in~\eqref{e:bdry-LA}-\eqref{e:bdry-LA-bis}-\eqref{e:bdry-LA-ter}, we have obtained 
\bna
L_{A}(h_n)=\frac{1}{i}a_{1}(0)|h_{n}\psi_{n}'(0^{+})|^{2}-\frac{1}{i}a_{1}(L)|h_{n}\psi_{n}'(L^{-})|^{2}+\grando{h_{n}^{1/2}}   \underset{n\to+\infty}{\to}  \frac1{i} ( a_1(0) \ell_0- a_1(L) \ell_L) ,
\ena
where we used~\eqref{limtra} in the limit. This concludes the proof of the lemma when combined with~\eqref{e:LA-on-one-side}.
  \enp
  We have used the following Lemma which is a 1D simpler version of \cite[Lemma 3.8]{GL:93}, and which proof relies on the elementary lemmata~\ref{lmtrace} and~\ref{lmtraceloin} below which sometimes use the specific properties (parity) for $\chi$ in~\eqref{e:def-chi}.
  \begin{lemma}  
    \label{lmtracejump}
Let $f \in L^2_{\comp}(\R)$ be such that, in $\mathcal{D}'(\R)$, we have
$$
f' = F + \alpha \delta_0 + \beta \delta_L ,  \quad \text{with}\quad  F\in  L^2(\R), \quad\alpha,\beta \in \C .
$$  
Then, with $\chi$ as in~\eqref{e:def-chi}, we have $f|_{(-\infty,0)} \in C^0([-\infty,0])$, $f|_{(0,L)} \in C^0([0,L])$, $f|_{(L,\infty)} \in C^0([L,\infty])$, together with 
$$
\big(\chi(\e D_{x})f\big)(0)=\frac{f(0^{+})+f(0^{-})}{2}+r , \quad \text{with} \quad |r| \leq C \e^{1/2} \left( \nor{F}{L^{2}(\R)}+ \nor{f}{L^{2}(\R)}\right) .
$$
  \end{lemma}
  \bnp
  The fact that $f$ is piecewise continuous follows from the fact that $f'$ is a Radon measure. 
  Using Lemma~\ref{lmtraceloin} below and a partition of the unity, we are reduced to the case where $f$ is supported in $(-L/2,L/2)$ and $\beta=0$.
  
  We define $g(x)=\frac{f(x)+f(-x)}{2}$. Then, $g$ is $C^{0}_{c}(\R)$  with $g'\in L^{2}(\R)$ and $\nor{g'}{L^{2}(\R)}\leq \nor{F}{L^{2}(\R)}$. We have $g(0)=\frac{f(0^{+})+f(0^{-})}{2}$. Using that $\chi$ is even and writing $\chi_{\e}(s)=\chi(\e s)$, we also have (denoting by $\check{\chi_{\e}}$ the inverse Fourier transform of $\chi_{\e}$)
  \bna
 \left[\chi(\e D_{x})g\right](0)=\frac{1}{2}\int_{\R_{y}}\left[f(-y)+f(y)\right]\check{\chi_{\e}}(y)dy=\frac{1}{2}\int_{\R_{y}}\left[\check{\chi_{\e}}(y)+\check{\chi_{\e}}(-y)\right]f(y)dy=  \left[\chi(\e D_{x})f\right](0).
  \ena
  We can conclude by applying Lemma~\ref{lmtrace} below to $g$.
  \enp
  \begin{lemma}
  \label{lmtrace}
  There is $C>0$ such that for all $f \in C^{0}_{c}(\R)$ with $f'\in L^{2}(\R)$, we have 
  $$\left|\big(\chi(\e D_{x})f\big)(0)-f(0)\right|\leq C \e^{1/2} \nor{f'}{L^{2}(\R)} , \quad \text{ for all }\eps>0 .
  $$
   \end{lemma}
\bnp
Denoting $\check{\chi_{\e}}$ the inverse Fourier transform of $\chi_{\e}$, we have $\check{\chi_{\e}}(\xi)=\frac{1}{\e}\check{\chi}(\e^{-1}\xi)$.
Since $\chi(0)=\chi_{\e}(0)=1$, we have $\int_{\R} \check{\chi_{\e}}(\xi) d\xi=1$ so that  
\begin{align*}
\big(\chi_{\e}(D_{x})f\big)(0)-f(0)&=\int_{\R_{y}}\left[f(-y)-f(0)\right]\check{\chi_{\e}}(y)dy=-\int_{\R_{y}}y\check{\chi_{\e}}(y)\int_{0}^{1}f'(-ty)dtdy \\
& =-\e \int_{\R_{x}}x\check{\chi}(x)\int_{0}^{1}f'(-t\e x)dtdx.
\end{align*}
We have by Cauchy-Schwarz inequality
\begin{align*}
\left|\big(\chi_{\e}(D_{x})f\big)(0)-f(0)\right|&\leq  \e \int_{0}^{1} \int_{\R_{x}}|x\check{\chi}(x)|\left|f'(-t\e x)\right|dtdx\leq \e \nor{x\check{\chi}}{L^{2}}  \int_{0}^{1} \left(\int_{\R_{x}}\left|f'(-t\e x)\right|^{2}dx\right)^{1/2}dt \\
& \leq C \e^{1/2}\int_{0}^{1}t^{-1/2} \left(\int_{\R_{s}}\left|f'(s)\right|^{2}ds\right)^{1/2}dt\leq C \e^{1/2} \nor{f'}{L^{2}(\R)} ,
\end{align*}
which is the sought estimate.
\enp
\begin{lemma}
\label{lmtraceloin} Let $c>0$ and $N \in \R_+$, then, there exists $C_{N,c}>0$ so that for all $f\in L^{2}(\R)$ so that $f=0$ a.e. in $(-c,c)$, we have
\bna
\left|\big(\chi(\e D_{x})f\big)(0)\right|\leq C_{N,c}\e^{N}\nor{f}{L^{2}} .
\ena  
\end{lemma}
\bnp
With the same notations as the proof of Lemma~\ref{lmtrace}, we have
\bna
\left|\big(\chi(\e D_{x})f\big)(0)\right|=\left|\int_{\R}f(-y)\check{\chi_{\e}}(y)dy\right|\leq \nor{y^{-N}f}{L^{2}}\nor{y^{N}\check{\chi_{\e}}(y)}{L^{2}}\leq \e^{N-1/2}c^{-N}\nor{f}{L^{2}}\nor{x^{N}\check{\chi}(x)}{L^{2}},
\ena
whence the sought result after having changed the value of $N$.
\enp

We have also used the following lemma to handle ``low-regularity'' potentials.
\begin{lemma}
\label{l:comm-pas-reg}
Assume $V \in C^0(\R)$ such that $V' \in L^\infty(\R)$ and $\chi \in C^\infty_c(\R)$. Then we have 
$$
[\chi(\eps D) , V(x)] \in \L(L^2(\R)) , \quad \text{ with } \quad \nor{[\chi(\eps D) , V(x)]}{\L(L^2)} \leq C_\chi \eps \nor{V'}{L^\infty(\R)}.
$$ 
\end{lemma}
\bnp
The operator $\chi(\eps D)$ is the convolution by $\frac{1}{\eps}\check{\chi}\left(\frac{\cdot}{\eps}\right)$ where $\check{\chi}$ is the inverse Fourier transform of $\chi$. Its kernel is therefore $\frac{1}{\e}\check{\chi}\left(\frac{x-y}{\e}\right)$ and the kernel of $[\chi(\eps D) , V(x)]$ is therefore $K_{\eps}(x,y)=\frac{1}{\e}\check{\chi}\left(\frac{x-y}{\e}\right)\left(V(y)-V(x)\right)$. The Schur Lemma and symmetry of the kernel in $(x,y)$ give
\begin{align*}
\nor{[\chi(\eps D) , V(x)]}{\L(L^2(\R))}&\leq \max\left[\sup_{x\in \R}\nor{K_{\eps}(x,y)}{L^{1}(\R_{y})},\sup_{y\in \R}\nor{K_{\eps}(x,y)}{L^{1}(\R_{x})}\right]\\
&\leq \frac{1}{\e} \sup_{x\in \R}\int_{\R_{y}}\left|\check{\chi}\left(\frac{x-y}{\e}\right)\right|\left|V(y)-V(x)\right|~dy \leq  \sup_{s\in \R}\int_{\R_{t}}\left|\check{\chi}\left(s-t\right)\right|\left|V(\eps t)-V(\eps s)\right|\\
& \leq \eps \nor{V'}{L^\infty(\R)}\sup_{s\in \R} \int_{\R_{t}}\left|\check{\chi}\left(s-t\right)(s-t)\right|~dt.
\end{align*}
This yields the expected result with $C_\chi=\nor{t\check{\chi}(t) }{L^{1}_{t}}$.
\enp
\subsection{Invariance properties near the boundary}
\label{s:invar-boundar}
Now, we will state the propagation at the boundary. We only consider the boundary problem at $x=0$, the problem near $x=L$ being handled similarly.
The following is a 1D version of \cite[Theorem~2.3]{GL:93}.
 \begin{lemma}
 \label{l:lemma-bdry-3-cas}
Under the assumptions of Proposition \ref{propmeasureapp}, with $\|r_n\|_{L^2} = \petito{1}$, we have
 \begin{itemize}
 \item (Elliptic case) if $\V(0)>0$: then $\ell_{0}=0$ and $\mu=0$ for $x$ close to $0$
\item (Glancing case) if $\V(0)=0$: then $H_{p}\mu=-\ell_{0}\delta_{x=0}\otimes \delta'_{\xi=0} $ for $x$ close to $0$
\item (Hyperbolic case) if $\V(0)<0$: then $H_{p}\mu=\frac{\ell_{0}}{2\sqrt{-\V(0)}}\delta_{x=0}\otimes (\delta_{\xi=\sqrt{-\V(0)}}- \delta_{\xi=-\sqrt{-\V(0)}}) $ for $x$ close to $0$. 
\end{itemize}
 \end{lemma} 
 Note that the simple 1D setting here allows to avoid the use of the Malgrange preparation theorem and provides with a self-contained elementary proof (as compared to~\cite[Theorem~2.3]{GL:93}).
 \begin{remark}
 Note that one recovers the equation in the glancing case $\V(0)=0$ by taking the limit $\sqrt{-\V(0)}\to 0$ in the equation obtained in the hyperbolic case.
  \end{remark}
   \bnp
   In the first elliptic case, that $\mu =0$ near $x=0$ is a consequence of Lemma \ref{lmsuppmu} together with $\{p=0\} \cap( \{0\}\times \R) =\emptyset$ if $\V(0)>0$. Applying Lemma \ref{lmderivbord} with $a_0=0$ and $a_{1}$ satisfying $a_{1}(0)=1$ and $\supp(a_1) \times \R \cap \supp(\mu) = \emptyset$ yields $\ell_0=0$.
   
   For the glancing case, we use Lemma \ref{lmsupHpmu} together with $\{\xi^{2}+\V(x)=0\}\cap \{0\}\times \R_{\xi}=\{(0,0)\}$. Classical Distribution theory implies that close to $x=0$, 
   \bna
   H_{p}\mu= q_{(0,0)}\delta_{(0,0)}+q_{(1,0)}\d_{x}\delta_{(0,0)}+q_{(0,1)}\d_{\xi}\delta_{(0,0)} ,
   \ena
   where $q_{\alpha}\in \C$.
Lemma \ref{lmderivbord} gives, for every $a(x,\xi)=a_{0}(x)+\xi a_{1}(x)$
   \bna
  \ell_{0}a_{1}(0)=-\left<\mu,H_{p}a\right>=q_{(0,0)}a_{0}(0)-q_{(1,0)}a_{0}'(0)-q_{(0,1)}a_{1}(0).
   \ena
   Since $a_{0}$ and $a_{1}$ are arbitrary smooth functions, we obtain $q_{(0,0)}=q_{(1,0)}=0$ and $q_{(0,1)}=-\ell_{0}$, so that $H_{p}\mu= -\ell_{0}\delta_{x=0}\otimes \delta'_{\xi=0}$, which is the sought result.
  
   For the hyperbolic case, Lemma \ref{lmsupHpmu} together with $$\{\xi^{2}+\V(x)=0\}\cap \{0\}\times \R_{\xi}=\{(0_{x},\sqrt{-\V(0)})\}\cup\{(0_{x},-\sqrt{-\V(0)})\}$$ 
   imply again that, close to $x=0$, 
      \begin{align}
      \label{formuleHpdirac}
   H_{p}\mu&= q^{+}_{(0,0)}\delta_{(0,\sqrt{-\V(0)})}+q^{+}_{(1,0)}\d_{x}\delta_{(0,\sqrt{-\V(0)})}+q^{+}_{(0,1)}\d_{\xi}\delta_{(0,\sqrt{-\V(0)})} \nonumber\\
   &\quad  +q^{-}_{(0,0)}\delta_{(0,-\sqrt{-\V(0)})}+q^{-}_{(1,0)}\d_{x}\delta_{(0,-\sqrt{-\V(0)})}+q^{-}_{(0,1)}\d_{\xi}\delta_{(0,-\sqrt{-\V(0)})}.
   \end{align}
   This time, Lemma \ref{lmderivbord} gives for every $a(x,\xi)=a_{0}(x)$
   \bna
0=-\left<\mu,H_{p}a\right>= \left<H_{p}\mu,a\right> =(q_{(0,0)}^{+}+q_{(0,0)}^{-})a_{0}(0)-(q_{(1,0)}^{+}+q_{(1,0)}^{-})a_{0}'(0).
   \ena
   This implies
   \begin{align}
   \label{e:qq}
   q_{(0,0)} := q_{(0,0)}^{+}=- q_{(0,0)}^- \quad \text{ and } \quad q_{(1,0)} := q_{(1,0)}^{+}=-q_{(1,0)}^-.
   \end{align} Then, Lemma \ref{lmderivbord} gives for every $a(x,\xi)=a_{1}(x)\xi$
      \bna
  \ell_{0}a_{1}(0)=-\left<\mu,H_{p}a\right>=\sqrt{-\V(0)}\left[2 a_{1}(0)q_{(0,0)} -2 a_{1}'(0)q_{(1,0)}\right]- (q^{+}_{(0,1)}+q^{-}_{(0,1)})a_{1}(0) .
   \ena
   This gives $q_{(1,0)}=0$ and 
   \begin{align}
   \label{e:ell-0-inter}
   \ell_{0}=\sqrt{-\V(0)}2 q_{(0,0)}-(q^{+}_{(0,1)}+q^{-}_{(0,1)}) .
   \end{align}
   To finish, we now choose $a(x,\xi)=p(x,\xi)b(x,\xi)$ as a test function, for $b\in C^\infty_c(\R^2)$, and obtain 
   \bna
   \left<\mu,H_{p}a\right>=  \left<\mu,p H_{p}b\right>+  \left<b\mu,H_{p}p\right>=0 ,
   \ena
   where we have used Lemma \ref{lmsuppmu} for the first term and $H_{p}p=0$ for the second. Applying again ~\eqref{formuleHpdirac} to this function $a$ and using the information we already have on the coefficients in~\eqref{formuleHpdirac}, we obtain using $a(0,\pm\sqrt{-\V(0)})=0$ (recall that $p=\xi^2+\V(x)$) that 
     $$
0=-\left<\mu,H_{p}a\right>= -q^{+}_{(0,1)}(\d_{\xi}a)(0,\sqrt{-\V(0)}) +q^{-}_{(0,1)}(\d_{\xi}a)(0,-\sqrt{-\V(0)}).
   $$
   But now for $a(x,\xi)=p(x,\xi)b(x,\xi)$, on the set $p=0$, we have 
   $$(\d_{\xi}a)(x,\xi)=(\d_{\xi} p)(x,\xi)b(x,\xi)+(\d_{\xi} b)(x,\xi)p(x,\xi)=2\xi b(x,\xi).
   $$ So, we deduce
   \bna
   0= -\sqrt{-\V(0)} q^{+}_{(0,1)}b (0,\sqrt{-\V(0)})-\sqrt{-\V(0)}q^{-}_{(0,1)}b(0,-\sqrt{-\V(0)}) .
   \ena
   Since $b$ is arbitrary and $\sqrt{-\V(0)}\neq 0$, we obtain $q^{+}_{(0,1)}=q^{-}_{(0,1)}=0$.
   This, together with~\eqref{e:ell-0-inter} implies that $\ell_{0}=2\sqrt{-\V(0)} q_{(0,0)}$ which, combined with~\eqref{formuleHpdirac},~\eqref{e:qq} and $q_{(1,0)}=0$, gives the expected result in the hyperbolic case.
      \enp

We now specify to the glancing and diffractive case at $x=0$.
\begin{lemma}
 If $\V(0)=0$ and $\V'(0)\leq 0$, then $\ell_{0}=0$. If moreover $\V'(0)<0$, then $\mu(\{(0,0)\})=0$.
\end{lemma}
\bnp
For this, we follow~\cite{BG:97}. We take $\chi \in C^\infty_c(-1,1)$ with $\chi=1$ in a neighborhood of $0$, $\chi\geq 0$ and $\int_\R \chi = 1$. Define $\tilde\chi(s) = \int_{-\infty}^s\chi \in C^\infty(\R)$ and test the identity~$H_{p}\mu=-\ell_{0}\delta_{x=0}\otimes \delta'_{\xi=0}$ obtained in Lemma~\ref{l:lemma-bdry-3-cas} with the function $a(x,\xi) = \chi(x/\alpha)\tilde\chi(\xi/\beta) \in C^\infty(\R^2)$ for $\alpha,\beta>0$. This yields (for $\alpha$ sufficiently small)
 \begin{align*}
\langle \mu , -\frac{2 \xi}{\alpha} \chi'(x/\alpha) \tilde\chi(\xi/\beta)\rangle + \langle \mu , \frac{\V'(x)}{\beta} \chi(x/\alpha) \tilde\chi'(\xi/\beta)\rangle = \frac{\ell_0}{\beta} \chi(0)\chi(0) =\frac{\ell_0}{\beta} .
\end{align*}
Multiplying by $\beta$, choosing $\alpha=\sqrt{\beta}$, and using dominated convergence yields, in the limit $\beta \to 0^+$ 
 \begin{align*}
  \grando{\sqrt{\beta}} + \langle \mu , \V'(x)  \chi(x/\sqrt{\beta}) \chi(\xi/\beta)\rangle = \ell_0  .
\end{align*}
Now taking the limit $\beta \to 0^+$ and using again dominated convergence implies $\V'(0) \mu (\{(0,0)\}) = \ell_0$. That $\V'(0)\leq 0$, $\mu\geq 0$ and $\ell_0 \geq 0$ implies that $\ell_0 = 0$. If moreover $\V'(0)<0$, then, we obtain $\mu(\{(0,0)\})=0$, which concludes the proof of the lemma.
     \enp 
      
\small
\bibliographystyle{alpha}
\bibliography{bibli}
\end{document}